\documentclass[aos,MSNbibl,nameyear,rotating,seceqn,dvips]{arximspdf}
\usepackage{dcolumn}
\usepackage{graphicx}

\doi{10.1214/13-AOS1099} 
\volume{41}
\issue{3}
\pubyear{2013}
\firstpage{1166}
\lastpage{1203}

\makeatletter
\newcolumntype{d}[1]{D{.}{.}{#1}}

\newcommand{\rrVert}{\Vert}
\newcommand{\rrvert}{\vert}
\newcommand{\llVert}{\Vert}
\newcommand{\llvert}{\vert}
\newcommand{\iint}{\int\!\!\!\int}
\newproclaim{ass}{Assumption}
\newtheorem{theorem}{Theorem}
\newtheorem{lemma}{Lemma}
\newtheorem{proposition}{Proposition}
\newproclaim{dgp}{DGP}
\makeatother

\begin{document}
\begin{frontmatter}

\title{A loss function approach to model specification testing and its relative efficiency}
\runtitle{Loss function approach}

\begin{aug}
\author[A]{\fnms{Yongmiao} \snm{Hong}\ead[label=e1]{yh20@cornell.edu}}
\and
\author[B]{\fnms{Yoon-Jin} \snm{Lee}\corref{}\ead[label=e2]{lee243@indiana.edu}}
\runauthor{Y. Hong and Y.-J. Lee}
\affiliation{Cornell University and Xiamen University, and Indiana University}
\address[A]{Department of Economics\\
\quad and Department of Statistical Science\\
Uris Hall\\
Cornell University\\
Ithaca, New York 14850\\
USA\\
and\\
Wang Yanan Institute for Studies in Economics\\
and MOE Key Laboratory of Economics\\
Xiamen University\\
Xiamen 361005, Fujian\\
China\\
\printead{e1}} 
\address[B]{Department of Economics\\
Indiana University\\
Wylie Hall\\
Bloomington, Indiana 47405\\
USA\\
\printead{e2}}
\end{aug}

\received{\smonth{1} \syear{2012}}
\revised{\smonth{10} \syear{2012}}

%
\begin{abstract}
The generalized likelihood ratio (GLR) test proposed by Fan, Zhang and
Zhang [\textit{Ann. Statist.}
\textbf{29} (2001) 153--193] and Fan and Yao [\textit{Nonlinear Time
Series: Nonparametric and Parametric Methods}
(2003) Springer] is a generally applicable nonparametric
inference procedure. In this paper, we show that although it inherits many
advantages of the parametric maximum likelihood ratio (LR) test, the GLR
test does not have the optimal power property. We propose a generally
applicable test based on loss functions, which measure discrepancies between
the null and nonparametric alternative models and are more relevant to
decision-making under uncertainty. The new test is asymptotically more
powerful than the GLR test in terms of Pitman's efficiency criterion. This
efficiency gain holds no matter what smoothing parameter and kernel function
are used and even when the true likelihood function is available for
the GLR
test.
\end{abstract}

%

\end{frontmatter}

\section{Introduction}\label{sec1}
The likelihood ratio (LR) principle is a generally applicable approach
to parametric hypothesis testing [e.g., Vuong (\citeyear{Vuo89})]. The maximum LR
test compares the best explanation of data under the alternative with
the best explanation under the null hypothesis. It is well known from
the Neyman--Pearson lemma that the maximum LR test has asymptotically
optimal power. Moreover, the LR statistic follows an asymptotic null
$\chi^{2}$ distribution with a known number of degrees of freedom,
enjoying the so-called Wilks phenomena that its asymptotic distribution
is free of nuisance parameters.

In parametric hypothesis testing, however, it is implicitly assumed
that the
family of alternative likelihood models contains the true model. When this
is not the case, one may fail to reject the null hypothesis
erroneously. In
many testing problems in practice, while the null hypothesis is well
formulated, the alternative is vague. Over the last two decades or so, there
has been a growing interest in nonparametric inference, namely, inference
for hypotheses on parametric, semiparametric and nonparametric models
against a nonparametric alternative. The nonparametric alternative is very
useful when there is no prior information about the true model. Because the
nonparametric alternative contains the true model at least for large
samples, it ensures the consistency of a test. Nevertheless, there have been
few generally applicable nonparametric inference principles. One naive
extension would be to develop a nonparametric maximum LR test similar
to the
parametric maximum LR test. However, the nonparametric maximum likelihood
estimator (MLE) usually does not exist, due to the well-known infinite
dimensional parameter problem [\citet{Bah58}, Le Cam (\citeyear{le90})]. Even if it
exists, it may be difficult to compute, and the resulting nonparametric
maximum LR test is not asymptotically optimal. This is because the
nonparametric MLE chooses the smoothing parameter automatically, which
limits the choice of the smoothing parameter and renders it impossible for
the test to be optimal.

\citet{FanZhaZha01} and \citet{FanYao03} proposed a generalized
likelihood ratio (GLR) test by replacing the nonparametric MLE with a
reasonable nonparametric estimator, attenuating the difficulty of the
nonparametric maximum LR test and enhancing the flexibility of the
test by
allowing for a range of smoothing parameters. The GLR test maintains the
intuitive feature of the parametric LR test because it is based on the
likelihoods of generating the observed sample under the null and alternative
hypotheses. It is generally applicable to various hypotheses involving a
parametric, semiparametric or nonparametric null model against a
nonparametric alternative. By a proper choice of the smoothing parameter,
the GLR test can achieve the asymptotically optimal rate of
convergence in
the sense of Ingster (\citeyear{Ing93N1,Ing93N2,Ing93N3}) and \citet{LepSpo99}.
Moreover, it enjoys the appealing Wilks phenomena that its asymptotic null
distribution is free of nuisance parameters and nuisance functions.

The GLR test is a nonparametric inference procedure based on the empirical
Kullback--Leibler information criterion (KLIC) between the null model
and a
nonparametric alternative model. This measure can capture any discrepancy
between the null and alternative models, ensuring the consistency of
the GLR
test. As Fan, Zhang and Zhang (\citeyear{FanZhaZha01}) and \citet{FanJia07} point out, it
holds an advantage over many discrepancy measures such as the $L_{2}$
and $L_{\infty}$ measures commonly used in the literature because for the latter
the choices of measures and weight functions are often arbitrary, and the
null distributions of the test statistics are unknown and generally depend
on nuisance parameters. We note that \citet{Rob91} developed a
nonparametric KLIC test for serial independence and White [(\citeyear
{Whi82}), page~17] also
suggested a nonparametric KLIC test for parametric likelihood models.

The GLR test assumes that stochastic errors follows some parametric
distribution which need not contain the true distribution. It is essentially
a nonparametric pseudo LR test. Azzalini, Bowman and H\"{a}rdle (\citeyear{AzzBowHar89}),
\citet{AzzBow90} and Cai, Fan and Yao (\citeyear{CaiFanYao00}) also proposed a
nonparametric pseudo-LR test for the validity of parametric regression
models.

In this paper, we show that despite its general nature and appealing
features, the GLR test does not have the optimal power property of the
classical LR test. We first propose a generally applicable nonparametric
inference procedure based on loss functions and show that it is
asymptotically more powerful than the GLR test in terms of Pitman's
efficiency criterion. Loss functions are often used in estimation, model
selection and prediction [e.g., Zellner (\citeyear{Zel86}), \citet{Phi96},
\citet{Wei96},
\citet{ChrDie97}, \citet{GiaWhi06}], but not in
testing. A loss function compares the models under the null and alternative
hypotheses by specifying a penalty for the discrepancy between the two
models. The use of a loss function is often more relevant to decision-making
under uncertainty because one can choose a loss function to mimic the
objective of the decision maker. In inflation forecasting, for example,
central banks may have asymmetric preferences which affect their optimal
policies [\citet{PeeNob98}]. They may be more concerned with
underprediction than overprediction of inflation rates. In financial risk
management, regulators may be more concerned with the left-tailed
distribution of portfolio returns than the rest of the distribution. In
these circumstances, it is more appropriate to choose an asymmetric loss
function to validate an inflation rate model and an asset return
distribution model. The admissible class of loss functions for our approach
is large, including quadratic, truncated quadratic and asymmetric linex loss
functions [Varian (\citeyear{Va75}), \citet{Zel86}]. They do not require any
knowledge of
the true likelihood, do not involve any choice of weights, and enjoy the
Wilks phenomena that its asymptotic distribution is free of nuisance
parameters and nuisance functions. Most importantly, the loss function test
is asymptotically more powerful than the GLR test in terms of Pitman's
efficiency criterion, regardless of the choice of the smoothing parameter
and the kernel function. This efficiency gain holds even when the true
likelihood function is available for the GLR test. Interestingly, all
admissible loss functions are asymptotically equally efficient under a
general class of local alternatives.

The paper is planned as follows. Section~\ref{sec2} introduces the framework and the
GLR principle. Section~\ref{sec3} proposes a class of loss function-based tests. For
concreteness, we focus on specification testing for time series regression
models, although our approach is applicable to other nonparametric testing
problems. Section~\ref{sec4} derives the asymptotic distributions of the loss
function test and the GLR test. Section~\ref{sec5} compares their relative efficiency
under a class of local alternatives. In Section~\ref{sec6}, a simulation study
compares the performance between two competing tests in finite samples.
Section~\ref{sec7} concludes the paper. All mathematical proofs are collected in an
\hyperref[app]{Appendix} and supplementary material [\citet{HonLee}].

\section{Generalized likelihood ratio test}\label{sec2}
Maximum LR tests are a generally applicable and powerful inference
method for most parametric testing problems. However, the classical LR
principle implicitly assumes that the alternative model contains the
true data generating process (DGP). This is not always the case in
practice. To ensure that the alternative model contains the true DGP,
one can use a nonparametric alternative model.

Recognizing the fact that the nonparametric MLE may not exist and so cannot
be a generally applicable method, Fan, Zhang and Zhang (\citeyear{FanZhaZha01}) and \citet{FanYao03} proposed the GLR principle as a generally applicable method for
nonparametric inference. The idea is to compare a suitable nonparametric
estimator with a restricted estimator under the null hypothesis via a LR
statistic. Specifically, suppose one is interested in whether a parametric
likelihood model $f_{\theta}$ is correctly specified for the unknown
density $f$ of the DGP, where $\theta$ is a finite-dimensional parameter.
The null hypothesis of interest is
%
\begin{equation}
\mathbb{H}_{0}\dvtx f=f_{\theta_{0}}\qquad\mbox{for some }\theta
_{0}\in\Theta, \label{eq2.1}
\end{equation}
where $\Theta$ is a parameter space. The alternative hypothesis is
%
\begin{equation}
\mathbb{H}_{A}\dvtx f\neq f_{\theta}\qquad\mbox{for all }\theta
\in\Theta. \label{eq2.2}
\end{equation}

In testing $\mathbb{H}_{0}$ versus $\mathbb{H}_{A}$, a nonparametric model
for $f$ can be used as an alternative, as also suggested in White
[(\citeyear{Whi82}),
page 17]. Suppose the log-likelihood function of a random sample is
$\hat{l}%
(f,\eta)$, where $\eta$ is a nuisance parameter. Under $\mathbb{H}_{0}$,
one can obtain the MLE $(\hat{\theta}_{0},\hat{\eta}_{0})$ by
maximizing the
model likelihood $\hat{l}(f_{\theta},\eta)$. Under the alternative $%
\mathbb{H}_{A}$, given $\eta$, one can obtain a reasonable smoothed
nonparametric estimator $\hat{f}_{\eta}$ of $f$. The nuisance
parameter $%
\eta$ can then be estimated by the profile likelihood; that is, to
find $%
\eta$ to maximize $l(\hat{f}_{\eta},\eta)$. This gives the maximum
profile likelihood $l(\hat{f}_{\hat{\eta}},\hat{\eta})$. The GLR test
statistic is then defined as
%
\begin{equation}
\lambda_{n}=l(\hat{f}_{\hat{\eta}},\hat{\eta})-l(f_{\hat{\theta}_{0}},
\hat{%
\eta}_{0}). \label{eq2.3}
\end{equation}
This is the difference of the log-likelihoods of generating the observed
sample under the alternative and null models. A large value of $\lambda_{n}$
is evidence against $\mathbb{H}_{0}$ since the alternative family of
nonparametric models is far more likely to generate the observed data.

The GLR test does not require knowing the true likelihood. This is appealing
since nonparametric testing problems do not assume that the underlying
distribution is known. For example, in a regression setting one usually does
not know the error distribution. Here, one can estimate model
parameters by
using a quasi-likelihood function $q(f_{\theta},\eta)$. The resulting GLR
test statistic\vadjust{\goodbreak} is then defined as%
%
\begin{equation}
\lambda_{n}=q(\hat{f}_{\hat{\eta}},\hat{\eta})-q(f_{\hat{\theta}_{0},}
\hat{%
\eta}_{0}). \label{eq2.4}
\end{equation}

The GLR approach is also applicable to the cases with unknown nuisance
functions. This can arise (e.g.) when one is interested in testing
whether a
function has an additive form which itself is still nonparametric. In this
case, one can replace $f_{\hat{\theta}_{0}}$\vspace*{1pt} by a nonparametric estimator
under the null hypothesis of additivity. \citet{Rob91} considers
such a
case in testing serial independence.

As a generally applicable nonparametric inference procedure, the GLR principle has been used to test a variety of models, including univariate
regression models [Fan, Zhang and Zhang (\citeyear{FanZhaZha01})], functional coefficient
regression models [Cai, Fan and Yao (\citeyear{CaiFanYao00})], spectral density models [Fan
and Zhang (\citeyear{FanZha04N1})], varying-coefficient partly linear regression models
[\citet{FanHua05}], additive models [\citet{FanJia05}],
diffusion models [\citet{FanZha03}] and partly linear additive
models [\citet{FanYao03}].
Analogous to the classical LR test statistic which follows an asymptotic
null $\chi^{2}$ distribution with a known number of degrees of
freedom, the
asymptotic distribution of the GLR statistic $\lambda_{n}$ is also a
$\chi
^{2}$ with a known large number of degrees of freedom, in the sense that
\[
r\lambda_{n}\simeq\chi_{\mu_{n}}^{2}
\]
as a sequence of constants $\mu_{n}\rightarrow\infty$ and some
constant $%
r>0;$ namely,
\[
\frac{r\lambda_{n}-\mu_{n}}{\sqrt{2\mu_{n}}}\stackrel{d}
{\rightarrow}%
N(0,1),
\]
where $\mu_{n}$ and $r$ are free of nuisance parameters and nuisance
functions, although they may depend on the methods of nonparametric
estimation and smoothing parameters. Therefore, the asymptotic distribution
of $\lambda_{n}$ is free of nuisance parameters and nuisance functions.
One can use $\lambda_{n}$ to make inference based on the known distribution
of $N(\mu_{n},2\mu_{n})$ or $\chi_{\mu_{n}}^{2}$ in large samples.
Alternatively, one can simulate the null distribution of $\lambda_{n}$ by
setting nuisance parameters at any reasonable values, such as the MLE
$\hat{%
\eta}_{0}$ or the maximum profile likelihood estimator $\hat{\eta}$ in
(\ref{eq2.3}).

The GLR test is powerful under a class of contiguous local
alternatives,%
\[
\mathbb{H}_{an}\dvtx f=f_{\theta_{0}}+n^{-\gamma}g_{n},
\]
where $\gamma>0$ is a constant and $g_{n}$ is an unspecified sequence of
smooth functions in a large class of function space. It has been shown [Fan, Zhang and Zhang (\citeyear{FanZhaZha01})] that when a local linear smoother is used
to estimate $%
f$ and the bandwidth is of order $n^{-2/9}$, the GLR test can detect local
alternatives with rate $\gamma=4/9$, which is optimal according to Ingster
(\citeyear{Ing93N1,Ing93N2,Ing93N3}).\looseness=-1

\section{A loss function approach}\label{sec3}
In this paper, we will show that while the GLR test enjoys many
appealing features of the classical LR test,\vadjust{\goodbreak} it does not have the
optimal power property of the classical LR test. We will propose a
class of loss function-based tests and show that they are
asymptotically more powerful than the GLR test under a class of local
alternatives. Loss functions measure discrepancies between the null and
alternative models and are more relevant to decision making under
uncertainty, because the loss function can be chosen to mimic the
objective function of the decision maker. The admissible loss functions
include but are not restricted to quadratic, truncated quadratic and
asymmetric linex loss functions. Like the GLR test, our tests are
generally applicable to various nonparametric inference problems, do
not involve choosing any weight function and their null asymptotic
distributions do not depend on nuisance parameters and nuisance
functions.

For concreteness, we focus on specification testing for time series
regression models. Regression modeling is one of the most important
statistical problems, and has been exhaustively studied, particularly
in the
i.i.d. contexts [e.g., H\"{a}rdle and Mammen (\citeyear{HarMam93})]. Focusing on
testing regression models will provide deep insight into our approach and
allow us to provide primitive regularity conditions for formal results.
Extension to time series contexts also allows us to expand the scope of
applicability of our tests and the GLR test. We emphasize that our approach
is applicable to many other nonparametric test problems, such as testing
parametric density models.

Suppose $\{X_{t},Y_{t}\} \in\mathbb{R}^{p+1}$ is a stationary time series
with finite second moments, where $Y_{t}$ is a scalar, $p\in\mathbb
{N}$ is
the dimension of vector $X_{t}$ and $X_{t}$ may contain exogenous and/or
lagged dependent variables. Then we can write%
%
\begin{equation}
Y_{t}=g_{0}(X_{t})+\varepsilon
_{t}, \label{eq3.1}
\end{equation}
where $g_{0}(X_{t})=E(Y_{t}|X_{t})$ and $E(\varepsilon_{t}|X_{t})=0$. The
fact that $E(\varepsilon_{t}|X_{t})=0$ does not imply that $\{
\varepsilon
_{t}\}$ is a martingale difference sequence. In a time series context,
$%
\varepsilon_{t}$ is often assumed to be i.i.d. $(0,\sigma^{2})$
and independent of $X_{t}$ [e.g., \citet{GaoGij08}]. This implies $
E(\varepsilon_{t}|X_{t})=0$ but not vice versa, and so it is overly
restrictive from a practical point of view. For example, $\varepsilon_{t}$
may display volatility clustering [e.g., \citet{Eng82}],%
\[
\varepsilon_{t}=z_{t}\sqrt{h_{t}},
\]
where $h_{t}=\alpha_{0}+\alpha_{1}\varepsilon_{t-1}^{2}$. Here, we
have $%
E(\varepsilon_{t}|X_{t})=0$ but $\{ \varepsilon_{t}\}$ is not
i.i.d. We will allow such an important feature, which is an empirical
stylized fact for high-frequency financial time series.

In practice, a parametric model is often used to approximate the unknown
function $g_{0}(X_{t})$. We are interested in testing validity of a
parametric model $g(X_{t},\theta)$, where $g(\cdot,\cdot)$ has a known
functional form, and $\theta\in\Theta$ is an unknown finite dimensional
parameter. The null hypothesis is
\[
\mathbb{H}_{0}\dvtx\Pr\bigl[g_{0}(X_{t})=g(X_{t},
\theta_{0})\bigr]=1\qquad\mbox{for some }%
\theta_{0}
\in\Theta
\]
versus the alternative hypothesis%
\[
\mathbb{H}_{A}\dvtx\Pr\bigl[g_{0}(X_{t})
\neq g(X_{t},\theta)\bigr]<1\qquad\mbox{for all }%
\theta\in
\Theta.
\]
An important example is a linear time series model%
\[
g(X_{t},\theta)=X_{t}^{\prime}\theta.
\]
This is called linearity testing in the time series literature [Granger and
Ter\"{a}svirta (\citeyear{GraTer93})]. Under $\mathbb{H}_{A}$, there exists neglected
nonlinearity in the conditional mean. For discussion on testing
linearity in
a time series context, see Granger and Ter\"{a}svirta (\citeyear{GraTer93}), Lee, White and Granger
(\citeyear{LeeWhiGra93}), Hansen (\citeyear{ha99}), Hjellvik and Tj{\o}stheim (\citeyear{HjeTjs96}) and
\citet{HonLee05}.

Because there are many possibilities for departures from a specific
functional form, and practitioners usually have no information about the
true alternative, it is desirable to construct a test of $\mathbb{H}_{0}$
against a nonparametric alternative, which contains the true function $%
g_{0}(\cdot)$ and thus ensures the consistency of the test against~$\mathbb{%
H}_{A}$. For this reason, the GLR test is attractive.

Suppose we have a random sample $\{Y_{t},X_{t}\}_{t=1}^{n}$ of size
$n\in
\mathbb{N}$. Assuming that the error $\varepsilon_{t}$ is i.i.d. $%
N(0,\sigma^{2})$, we obtain the conditional quasi-log-likelihood function
of $Y_{t}$ given $X_{t}$ as follows:%
%
\begin{equation}
\hat{l}\bigl(g,\sigma^{2}\bigr)=-\frac{n}{2}\ln\bigl(2\pi
\sigma^{2}\bigr)-\frac{1}{2\sigma
^{2}}\sum_{t=1}^{n}
\bigl[Y_{t}-g(X_{t})\bigr]^{2}. \label{eq3.2}
\end{equation}
Let $\hat{g}(x)$ be a consistent local smoother for $g_{0}(x)$.
Examples of $%
\hat{g}(x)$ include the Nadaraya--Watson estimator [H\"{a}rdle (\citeyear{Har90}),
\citet{LiRac07}, \citet{PagUll99}] and local linear estimator
[\citet{FanYao03}]. Substituting $\hat{g}(X_{t})$ into (\ref{eq3.2}), one
obtains the likelihood
of generating the observed sample $\{Y_{t},X_{t}\}_{t=1}^{n}$ under
$\mathbb{%
H}_{A}$,%
%
\begin{equation}
\hat{l}\bigl(\hat{g},\sigma^{2}\bigr)=-\frac{n}{2}\ln\bigl(2
\pi\sigma^{2}\bigr)-\frac{1}{%
2\sigma^{2}}\operatorname{SSR}_{1}, \label{eq3.3}
\end{equation}
where $\operatorname{SSR}_{1}$ is the sum of squared residuals of the nonparametric model;
namely,
\[
\operatorname{SSR}_{1}=\sum_{t=1}^{n}
\bigl[Y_{t}-\hat{g}(X_{t})\bigr]^{2}.
\]

Maximizing the likelihood in (\ref{eq3.3}) with respect to nuisance parameter $%
\sigma^{2}$ yields $\hat{\sigma}^{2}=n^{-1}\operatorname{SSR}_{1}$. Substituting this
estimator in (\ref{eq3.3}) yields the following likelihood:%
%
\begin{equation}
\hat{l}\bigl(\hat{g},\hat{\sigma}^{2}\bigr)=-\frac{n}{2}\ln
(\operatorname{SSR}_{1})-\frac{n}{2} \bigl[ 1+\ln(2\pi/n) \bigr].
\label{eq3.4}
\end{equation}

Using a similar argument and maximizing the model quasi-likelihood function
with respect to $\theta$ and $\sigma^{2}$ simultaneously,\vadjust{\goodbreak} we can obtain
the parametric maximum quasi-likelihood under $\mathbb{H}_{0}$,%
%
\begin{equation}
\hat{l}\bigl(\hat{g}_{\hat{\theta}_{0}},\hat{\sigma}_{0}^{2}
\bigr)=-\frac{n}{2}\ln \operatorname{SSR}_{0}-\frac{n}{2}\ln\bigl[1+\ln
(2\pi/n)\bigr], \label{eq3.5}
\end{equation}
where $(\hat{\theta}_{0},\hat{\sigma}_{0}^{2})$ are the MLE under
$\mathbb{H}%
_{0}$, and $\operatorname{SSR}_{0}$ is the sum of squared residuals of the parametric
regression model, namely,
\[
\operatorname{SSR}_{0}=\sum_{t=1}^{n} \bigl[
Y_{t}-g_{0}(X_{t},\hat{\theta}_{0})
\bigr] ^{2}.
\]
Given the i.i.d. $N(0,\sigma^{2})$ assumption for $\varepsilon
_{t}, $ $\hat{\theta}_{0}$ is the least squares estimator that
minimizes $%
\operatorname{SSR}_{0}$.

Thus, the GLR statistic is defined as
%
\begin{equation}
\lambda_{n}=\hat{l}\bigl(\hat{g},\hat{\sigma}^{2}\bigr)-
\hat{l}\bigl(\hat{g}_{\hat{\theta}%
_{0}},\hat{\sigma}_{0}^{2}
\bigr)=\frac{n}{2}\ln(\operatorname{SSR}_{0}/\operatorname{SSR}_{1}). \label{eq3.6}
\end{equation}

Under the i.i.d. $N(0,\sigma^{2})$ assumption for $\varepsilon
_{t}$, $\lambda_{n}$ is asymptotically equivalent to the $F$ test statistic
%
\begin{equation}
F=\frac{\operatorname{SSR}_{0}-\operatorname{SSR}_{1}}{\operatorname{SSR}_{1}}. \label{eq3.7}
\end{equation}
The latter has been proposed by Azzalini, Bowman and
H{\"a}rdle (\citeyear{AzzBowHar89}), Azzalini and Bowman
(\citeyear{AzzBow93}), \citet{HonWhi95} and \citet{FanLi02} in
i.i.d. contexts. The asymptotic equivalence between the GLR and $F$
tests can be seen from a Taylor series expansion of $\lambda_{n}$,%
\[
\lambda_{n}=\frac{n}{2}\cdot F+\mbox{Remainder}.
\]

We now propose an alternative approach to testing $\mathbb{H}_{0}$
versus $%
\mathbb{H}_{A}$ by comparing the null and alternative models via a loss
function $D\dvtx\mathbb{R}^{2}\rightarrow\mathbb{R}$, which measures the
discrepancy between the fitted values $\hat{g}(X_{t})$ and $g(X_{t},\hat
{%
\theta}_{0})$,%
%
\begin{equation}
Q_{n}=\sum_{t=1}^{n}D \bigl[
\hat{g}(X_{t}),g_{0}(X_{t},\hat{
\theta}_{0})%
\bigr]. \label{eq3.8}
\end{equation}
Intuitively, the loss function gives a penalty whenever the parametric model
overestimates or underestimates the true model. The latter is consistently
estimated by a nonparametric method.

A specific class of loss functions $D(\cdot,\cdot)$ is given by $%
D(u,v)=d(u-v)$, where $d(z)$ has a unique minimum at $0$, and is
monotonically nondecreasing as $|z|$ increases. Suppose $d(\cdot)$ is twice
continuously differentiable at $0$ with $d(0)=0,d^{\prime}(0)=0$ and $
0<d^{\prime\prime}(0)<\infty$. The condition of $d^{\prime}(0)=0$
implies that the first-order term in the Taylor expansion of $d(\cdot)$
around 0 vanishes to $0$ identically. This class of loss functions
$d(\cdot
) $ has been called a generalized cost-of-error function in the literature
[e.g., Pesaran and Skouras (\citeyear{PeSk01}), Granger (\citeyear{Gr99}), Christoffersen and Diebold
(\citeyear{ChrDie97}), Granger and Pesaran (\citeyear{GrPe00}), \citet
{Wei96}]. The loss function is
closely related to decision-based evaluation, which assesses the economic
value of forecasts to a particular decision maker or group of decision
makers. For example, in risk management the extreme values of portfolio
returns are of particular interest to regulators, while in macroeconomic
management the values of inflation or output growth, in the middle of the
distribution, may be of concern to central banks. A suitable choice of loss
function can mimic the objective of the decision maker.

Infinitely many loss functions $d(\cdot)$ satisfy the aforementioned
conditions, although they may have quite different shapes. To
illustrate the
scope of this class of loss functions, we consider some examples. The first
example of $d(\cdot)$ is the popular quadratic loss function%
%
\begin{equation}
d(z)=z^{2}. \label{eq3.9}
\end{equation}
This delivers a statistic based on the sum of squared differences between
the fitted values of the null and alternative models,%
%
\begin{equation}
\hat{L}_{n}^{2}=\sum_{t=1}^{n}
\bigl[\hat{g}(X_{t})-g_{0}(X_{t},\hat{
\theta}%
_{0})\bigr]^{2}. \label{eq3.10}
\end{equation}
This statistic is used in \citet{HonWhi95} and \citet
{HorSpo01} in an i.i.d. setup. It is also closely related to the
statistics proposed by H\"{a}rdle and Mammen (\citeyear{HarMam93}) and Pan, Wang and Yao
(\citeyear{PanWanYao07}) but different from their statistics, $\hat{L}_{n}^{2}$ in (\ref{eq3.10}) does
not involve any weighting which suffers from the undesirable feature as
pointed out in \citet{FanJia07}.

A second example of $d(\cdot)$ is the truncated quadratic loss function
%
\begin{equation}
d(z)= \cases{
\frac{1}{2}z^{2}, & \quad$\mbox{if }|z|\leq
c,$
\vspace*{2pt}\cr
 c|z|-\frac{1}{2}c^{2}, & \quad$\mbox{if }|z|>c,$}
%
\end{equation}
where $c$ is a prespecified constant. This loss function is used in
robust $%
M $-estimation. It is expected to deliver a test robust to outliers
that may
cause extreme discrepancies between two estimators.

The quadratic and truncated quadratic loss functions give equal penalty to
overestimation and underestimation of same magnitude. They cannot capture
asymmetric loss features that may arise in practice. For example, central
banks may be more concerned with underprediction than overprediction of
inflation rates. For another example, in providing an estimate of the market
value of a property of the owner, a real estate agent's underestimation and
overestimation may have different consequences. If the valuation is in
preparation for a future sale, underestimation may lead to the owner losing
money and overestimation to market resistance [Varian (\citeyear{Va75})].

\begin{figure}

\includegraphics{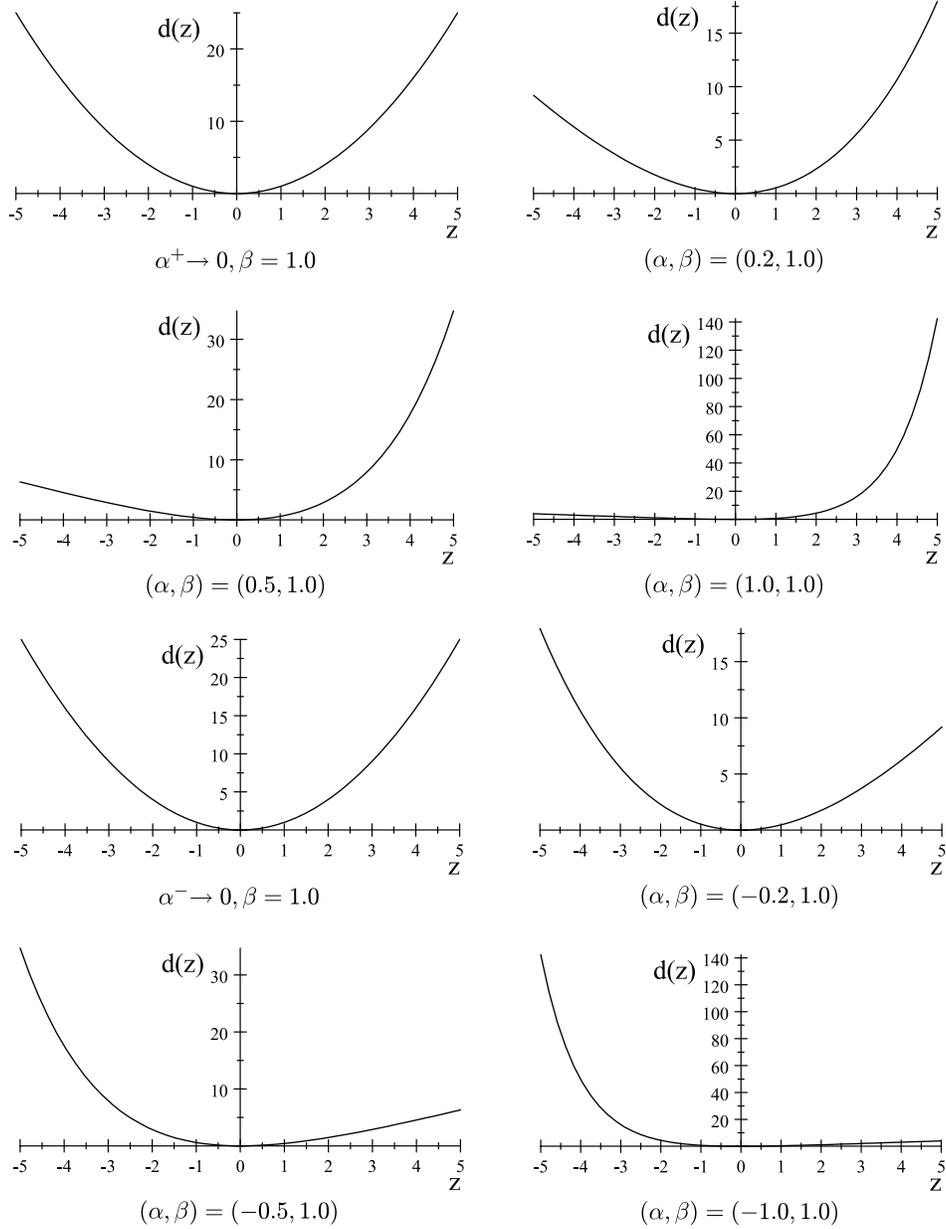}

  \caption{The LINEX loss function $d(z)=\frac{\beta}{\alpha^{2}}[\exp(\alpha
z)-(1+\alpha z)]$.}\label{fig1}
\end{figure}

The above examples motivate using an asymmetric loss function for model
validation. Examples of asymmetric loss functions are a class of so-called
linex functions%
%
\begin{equation}
d(z)=\frac{\beta}{\alpha^{2}}\bigl[\exp(\alpha z)-(1+\alpha z)\bigr
]. \label{eq3.12}
\end{equation}
For each pair of parameters $(\alpha,\beta)$, $d(z)$ is an asymmetric loss
function. Here, $\beta$ is a scale factor, and $\alpha$ is a shape
parameter. The magnitude of $\alpha$ controls the degree of asymmetry, and
the sign of $\alpha$ reflects the direction of asymmetry. When $\alpha<0$,
$d(z)$ increases almost exponentially if $z<0$, and almost linearly if $z>0$,
and conversely when $\alpha>0$. Thus, for this loss function,
underestimation is more costly than overestimation when $\alpha<0$,
and the
reverse is true when $\alpha>0$. For small values of $|\alpha|$,
$d(z)$ is
almost symmetric and not far from a quadratic loss function. Indeed, if
$%
\alpha\rightarrow0$, the linex loss function becomes a quadratic loss
function%
\[
d(z)\rightarrow\frac{\beta}{2}z^{2}.
\]
However, when $|\alpha|$ assumes appreciable values, the linex loss
function $d(z)$ will be quite different from a quadratic loss function.
Thus, the linex loss function can be viewed as a generalization of the
quadratic loss function allowing for asymmetry. This function was first
introduced by Varian (\citeyear{Va75}) for real estate assessment. \citet{Zel86}
employs it in the analysis of several central statistical estimation and
prediction problems in a Bayesian framework. \citet{GraPes} also
use it to evaluate density forecasts, and \citet{ChrDie97}
analyze the optimal prediction problem under this loss function. Figure~\ref{fig1}
shows the shapes of the linex function for a variety of choices of
$(\alpha,\beta)$.

Our loss function approach is by no means only applicable to regression
functions. For example, in such contexts as probability density\vspace*{1pt} and spectral
density estimation, one may compare two nonnegative density estimators,
say $%
f_{\hat{\theta}}$ and $\hat{f}$, using the Hellinger loss function%
%
\begin{equation}
D(f_{\hat{\theta}},\hat{f})= ( 1-\sqrt{f_{\hat{\theta}}/\hat{f}} )
^{2}. \label{eq3.13}
\end{equation}
This is expected to deliver a consistent robust test for $\mathbb
{H}_{0}$ of
(\ref{eq2.1}). Our approach covers this loss function as well, because when
$f_{\hat{%
\theta}}$ and $\hat{f}$ are close under $\mathbb{H}_{0}$ of~(\ref{eq2.1}), we have
\[
D(f_{\hat{\theta}},\hat{f})=\frac{1}{4} \biggl( \frac{f_{\hat{\theta
}}-\hat{f}%
}{\hat{f}}
\biggr) ^{2}+\mbox{Remainder},
\]
where the first-order term in the Taylor expansion vanishes to 0
identically. Interestingly, our approach does not apply to the KLIC loss
function%
%
\begin{equation}
D(f_{\hat{\theta}},\hat{f})=-\ln( f_{\hat{\theta}}/\hat{f} ), %
\label{eq3.14}
\end{equation}
which delivers the GLR $\lambda_{n}$ in (\ref{eq2.3}). This is because the Taylor
expansion of (\ref{eq3.14}) yields
%
\begin{equation}
D(f_{\hat{\theta}},\hat{f})=- \biggl( \frac{f_{\hat{\theta}}-\hat
{f}}{\hat{f}}%
\biggr) +
\frac{1}{2} \biggl( \frac{f_{\hat{\theta}}-\hat{f}}{\hat{f}} \biggr)
^{2}+\mbox{Remainder}, \label{eq3.15}
\end{equation}
where the first-order term in the Taylor expansion does not vanish to 0
identically. Hence, the first two terms in (\ref{eq3.15}) jointly determine the
asymptotic distribution of the GLR statistic $\lambda_{n}$. As will be seen
below, the presence of the first-order term in the Taylor expansion of the
KLIC loss function in (\ref{eq3.14}) leads to an efficiency loss compared to our
loss function approach for which the first-order term of a Taylor expansion
is identically 0 under the null.

\section{Asymptotic null distribution}\label{sec4}
Using a local fit with kernel $K\dvtx%
\mathbb{R}\rightarrow\mathbb{R}$ and bandwidth $h\equiv h(n)$, one could
obtain a nonparametric regression estimator\vspace*{1pt} $\hat{g}(\cdot)$ and
compare it
to the parametric model $g(\cdot,\hat{\theta}_{0})$ via a loss function,
where $\hat{\theta}_{0}$ is a consistent estimator for $\theta_{0}$
under $%
\mathbb{H}_{0}$. To avoid undersmoothing [i.e., to choose $h$ such that the
squared bias of $\hat{g}(\cdot)$ vanishes to 0 faster than the
variance of $%
\hat{g}(\cdot)$], we estimate the conditional mean of the estimated
parametric residual%
\[
\hat{\varepsilon}_{t}=Y_{t}-g(X_{t},\hat{
\theta}_{0}),
\]
and compare it to a zero function $E(\varepsilon_{t}|X_{t})=0$
(implied by $%
\mathbb{H}_{0})$ via a loss function criterion%
%
\begin{equation}
\hat{Q}_{n}=\sum_{t=1}^{n}D
\bigl[\hat{m}_{h}(X_{t}),0\bigr]=\sum
_{t=1}^{n}d\bigl[\hat{m}%
_{h}(X_{t})-0\bigr]=\sum_{t=1}^{n}d
\bigl[\hat{m}_{h}(X_{t})\bigr], \label{eq4.1}
\end{equation}
where $\hat{m}_{h}(X_{t})$ is a nonparametric estimator for
$E(\varepsilon
_{t}|X_{t})$. This is essentially a bias-reduction device. It is
proposed in
H\"{a}rdle and Mammen (\citeyear{HarMam93}) and also used in \citet{FanJia07} for the
GLR test. This device helps remove the bias of nonparametric estimation
because there is no bias under $\mathbb{H}_{0}$ when we estimate the
conditional mean of the estimated model residuals. We note that the
bias-reduction device does not lead to any efficiency gain of the loss
function test. The same efficiency gain of the loss function approach over
the GLR approach is obtained even when we compare estimators for $%
E(Y_{t}|X_{t})$. In the latter case, however, more restrictive
conditions on
the bandwidth $h$ are required to ensure that the bias vanishes sufficiently
fast under $\mathbb{H}_{0}$.

For simplicity, we use the Nadaraya--Watson estimator%
%
\begin{equation}
\hat{m}_{h}(x)=\frac{n^{-1}\sum_{t=1}^{n}\hat{\varepsilon}_{t}\mathbf{K}
_{h}(x-X_{t})}{n^{-1}\sum_{t=1}^{n}\mathbf{K}_{h}(x-X_{t})}, \label{eq4.2}
\end{equation}
where $X_{t}=(X_{1t},\ldots,X_{pt})^{\prime}$, $x=(x_{1},\ldots,x_{p})^{\prime}$, and%
\[
\mathbf{K}_{h}(x-X_{t})=h^{-p}\prod
_{i=1}^{p}K\bigl[h^{-1}(x_{i}-X_{it})
\bigr].
\]
We note that a local polynomial estimator could also be used, with the same
asymptotic results.

To derive the null limit distributions of the loss function test based
on $%
\hat{Q}_{n}$ in~(\ref{eq4.1}) and the GLR statistic $\lambda_{n}$ in a time series
context, we provide the following regularity conditions:\vspace*{-1pt}

\renewcommand{\theass}{A.\arabic{ass}}
\begin{ass}\label{assa1}
(i) For each $n\in\mathbb{N}$, $%
\{(Y_{t},X_{t}^{\prime})^{\prime}\in\mathbb{R}^{p+1}$, $t=1,\ldots,n\}$,
$p\in\mathbb{N}$, is a stationary and absolutely regular mixing process
with mixing coefficient $\beta(j)\leq C\rho^{j}$ for all $j\geq0$,
where $%
\rho\in(0,1)$, and $C\in(0,\infty);$ (ii)~$E|Y_{t}|^{8+\delta}<C$ for
some $\delta\in(0,\infty);$ (iii) $X_{t}$ has a compact support
$\mathbb{G%
}\subset\mathbb{R}^{p}$ with marginal probability density $C^{-1}\leq
f(x)\leq C$ for all $x$ in $\mathbb{G}$, and $f(\cdot)$ is twice
continuously differentiable on $\mathbb{G}$; (iv) the joint probability
density of $(X_{t},X_{t-j})$, $f_{j}(x,y)\leq C$ for all $j>0$ and all $
x,y\in\mathbb{G}$, where $C\in(0,\infty)$ does not depend on $j;$
(v) $%
E(X_{it}^{4(1+\eta)})\leq C$ for some $\eta\in(0,\infty)$, $1\leq
i\leq
p;$ (vi) $\operatorname{var}(\varepsilon_{t})=\sigma^{2}$ and $\sigma
^{2}(x)=E(\varepsilon_{t}^{2}|X_{t}=x)$ is continuous on $\mathbb{G}$.\vspace*{-1pt}
\end{ass}

\begin{ass}\label{assa2}
(i) For each $\theta\in\Theta$, $%
g(\cdot,\theta)$ is a measurable function of~$X_{t};$ (ii)~with
probability one, $g(X_{t},\cdot)$ is twice continuously differentiable with
respect to $\theta\in\Theta$, with $E\sup_{\theta\in\Theta
_{0}}\llVert \frac{\partial}{\partial\theta}g(X_{t},\theta
)\rrVert ^{4+\delta}\leq C$ and $E\sup_{\theta\in\Theta
_{0}}\hspace*{-1pt}\llVert \frac{\partial}{\partial\theta\, \partial\theta
^{\prime}%
}g(X_{t},\theta)\rrVert ^{4}\hspace*{-1pt}\leq C$, where $\Theta_{0}$ is a small
neighborhood of $\theta_{0}$ in~$\Theta$.\vspace*{-1pt}
\end{ass}

\begin{ass}\label{assa3}There exists a sequence of constants $%
\theta_{n}^{\ast}\in$ int$(\Theta) $ such that $n^{1/2}(\hat{\theta}
_{0}-\theta_{n}^{\ast})=O_{p}(1)$, where $\theta_{n}^{\ast}=\theta_{0}$
under $\mathbb{H}_{0}$ for all $n\geq1$.\vspace*{-1pt}
\end{ass}

\begin{ass}\label{assa4}
The kernel $K\dvtx\mathbb{R}\rightarrow
[0,1]$ is a prespecified bounded symmetric probability
density which satisfies the Lipschitz condition.\vspace*{-1pt}
\end{ass}

\begin{ass}\label{assa5}
$d\dvtx\mathbb{R}\rightarrow\mathbb{R}^{+}$
has a unique minimum at 0 and $d(z)$ is monotonically nondecreasing as $
|z|\rightarrow\infty$. Furthermore, $d(z)$ is twice continuously
differentiable at 0 with $d(0)=0,d^{\prime}(0)=0,D\equiv\frac{1}{2}%
d^{\prime\prime}(0)\in(0,\infty)$ and $|d^{\prime\prime
}(z)-d^{\prime
\prime}(0)|\leq C|z|$ for any $z$ near 0.\vspace*{-1pt}
\end{ass}

Assumptions~\ref{assa1} and~\ref{assa2} are conditions on the DGP. For each $t$, we
allow $%
(X_{t},Y_{t})$ to depend on the sample size $n$. This facilitates local
power analysis. For notational simplicity, we have suppressed the dependence
of $(X_{t},Y_{t})$ on $n$. We also allow time series data with weak serial
dependence. For the $\beta$-mixing condition, see, for example, Doukhan
(\citeyear{Dou94}). The
compact support for regressor $X_{t}$ is assumed in \citet
{FanZhaZha01} for the GLR test to avoid the awkward problem of tackling
the KLIC
function. This assumption allows us to focus on essentials while maintaining
a relatively simple treatment. It could be relaxed in several ways. For
example, we could impose a weight function $\mathbf{1}(|X_{t}|<C_{n})$ in
constructing $Q_{n}$ and $\lambda_{n}$, where $\mathbf{1}(\cdot)$ is the
indicator function, and $C_{n}$ can be either fixed or grow at a suitable
rate as the sample size $n\rightarrow\infty$.

Assumption~\ref{assa3} requires a $\sqrt{n}$-consistent estimator $\hat{\theta}_{0}$
under $\mathbb{H}_{0}$, which need not be asymptotically most
efficient. It
can be the conditional least squares or quasi-MLE.\vadjust{\goodbreak}\vspace*{1pt} Also, we do not need to
know the asymptotic expansion structure of $\hat{\theta}_{0}$ because the
sampling variation in $\hat{\theta}_{0}$ does not affect the limit
distribution of~$\hat{Q}_{n}$. We can estimate $\hat{\theta}_{0}$ and
proceed as if it were equal to $\theta_{0}$. The replacement of $\hat
{\theta%
}_{0}$ with $\theta_{0}$ has no impact on the limit distribution of
$\hat{Q}%
_{n}$.

We first derive the limit distribution of the loss function test
statistic.%

\begin{theorem}[(Loss function test)]\label{th1}
Suppose
Assumptions~\ref{assa1}--\ref{assa5} hold, $h\propto n^{-\omega}$ for
$\omega
\in(0,1/2p)$ and $p<4$. Define $q_{n}=\hat{Q}_{n}/\hat{\sigma
}%
_{n}^{2}$ where $\hat{Q}_{n}$ is given in (\ref{eq4.1}) and
$\hat{%
\sigma}_{n}^{2}=n^{-1}\operatorname{SSR}_{1}=n^{-1}\sum_{t=1}^{n}[\hat{\varepsilon
}_{t}-%
\hat{m}_{h}(X_{t})]^{2}$. Then \textup{(i)} under $\mathbb{H}_{0}$, $%
s(K)q_{n}\stackrel{d}{\simeq}\chi_{\nu_{n}}^{2}$as $n\rightarrow
\infty$, in the sense that
\[
\frac{s(K)q_{n}-\nu_{n}}{\sqrt{2\nu_{n}}}\stackrel{d}
{\longrightarrow}%
N(0,1),
\]
where $s(K)=\sigma^{2}a(K)\int\sigma^{2}(x)\,dx/[Db(K)\int
\sigma
^{4}(x)\,dx]$, $\nu_{n}=a^{2}(K)\times\break  [ \int\sigma^{2}(x)\,dx ]
^{2}/[h^{p}b(K)\int\sigma^{4}(x)\,dx]$, $a(K)=\int\mathbf
{K}^{2}(\mathbf{u}%
)\,d\mathbf{u}$, $b(K)=\int[\int\mathbf{K}(\mathbf
{u}+\break\mathbf{v})%
\mathbf{K}(\mathbf{v})\,d\mathbf{v]}^{2}\,d\mathbf{u}$, $\mathbf
{K}(\mathbf{u}%
)=\prod_{i=1}^{p}K(u_{i}),\mathbf{u}=(u_{1},\ldots,u_{p})^{\prime}$.

\textup{(ii)} Suppose in addition $\operatorname{var}(\varepsilon_{t}|X_{t})=\sigma
^{2}$ almost surely. Then $s(K)=a(K)/\break  [Db(K)]$ and $\nu
_{n}=\Omega a^{2}(K)/[h^{p}b(K)]$, where $\Omega$ is the
Lebesgue's measure of the support of $X_{t}$.
\end{theorem}

Theorem~\ref{th1} shows that under (and only under) conditional homoskedasticity,
the factors $s(K)$ and $\nu_{n}$ do not depend on nuisance parameters and
nuisance functions. In this case, the loss function test statistic $q_{n}$,
like the GLR statistic $\lambda_{n}$, also enjoys the Wilks phenomena that
its asymptotic distribution does not depend on nuisance parameters and
nuisance functions. This offers great convenience in implementing the loss
function test.

We note that the condition on the bandwidth $h$ is relatively mild. In
particular, no undersmoothing is required. This occurs because we estimate
the conditional mean of the residuals of the parametric model $%
g(X_{t},\theta)$. If we directly compared a nonparametric estimator of
$%
E(Y_{t}|X_{t})$ with $g(X_{t},\theta)$, we could obtain the same asymptotic
distribution for $q_{n}$, but under a more restrictive condition on $h$ in
order to remove the effect of the bias. For simplicity, we consider the case
with $p<4$. A higher dimension $p$ for $X_{t}$ could be allowed by suitably
modifying factors $s(K)$ and $\nu_{n}$, but with more tedious expressions.

Theorem~\ref{th1} also holds for the statistic $q_{n}^{0}=\hat{Q}_{n}/\hat
{\sigma}%
_{n,0}^{2}$, where $\hat{\sigma}_{n,0}^{2}=n^{-1}\operatorname{SSR}_{0}$, which is expected
to have better sizes than $q_{n}$ in finite samples under $\mathbb{H}_{0}$
when using asymptotic theory. However, $q_{n}$ may have better power
than $%
q_{n}^{0}$ because $\operatorname{SSR}_{0}$ may be substantially larger than $\operatorname{SSR}_{1}$
under $\mathbb{H}_{A}$.

To compare the $q_{n}$ and GLR tests, we have to derive the asymptotic
distribution of the GLR statistic $\lambda_{n}$ in a time series context,
a formal result not available in the previous literature, although the GLR
test has been widely applied in the time series context [\citet
{FanYao03}].\vadjust{\goodbreak}

\begin{theorem}[(GLR test in time series)]\label{th2}
Suppose
Assumptions~\ref{assa1}--\ref{assa5} hold, $p<4$, and $h\propto
n^{-\omega}$ for $\omega\in(0,1/2p)$, and $p<4$.
Define $\lambda_{n}$ as in (\ref{eq3.6}), where $\operatorname{SSR}_{1}=\sum_{t=1}^{n}[\hat{\varepsilon}_{t}-\hat{m}_{h}(X_{t})]^{2}$,
$\operatorname{SSR}_{0}=\sum_{t=1}^{n}\hat{\varepsilon}_{t}^{2}$ and
$\hat{\varepsilon}_{t}=Y_{t}-m(X_{t},\hat{\theta})$. Then \textup{(i)}~under
$\mathbb{H}_{0}$, $r(K)\lambda_{n}\stackrel{d}{\simeq}\chi_{\mu_{n}}^{2}$
as $n\rightarrow\infty$, in the sense that
\[
\frac{r(K)\lambda_{n}-\mu_{n}}{\sqrt{2\mu_{n}}}\stackrel{d} {%
\longrightarrow}N(0,1),
\]
where $r(K)=\sigma^{2}c(K)\int\sigma^{2}(x)\,dx/[d(K)\int
\sigma
^{4}(x)\,dx]$, $\mu_{n}= [ c(K)\int\sigma^{2}(x)\,dx ]
^{2}/\break [h^{p}\,d(K)\int\sigma^{4}(x)\,dx]$, $c(K)=\mathbf{K}(0)-\frac
{1}{2}\int
\mathbf{K}^{2}(\mathbf{u})\,d\mathbf{u}$, $d(K)=\int[\mathbf
{K}(%
\mathbf{u})-\break\frac{1}{2}\int\mathbf{K}(\mathbf{u}+\mathbf{v})\mathbf{K}(
\mathbf{v})\,d\mathbf{v}]^{2}\,d\mathbf{u}$, $\mathbf{K}(\mathbf{u}%
)=\prod_{i=1}^{p}K(u_{i})$, $\mathbf{u}=(u_{1},\ldots,u_{p})^{\prime}$.

\textup{(ii)} Suppose in addition $\operatorname{var}(\varepsilon_{t}|X_{t})=\sigma
^{2}$
almost surely. Then $r(K)=c(K)/\break d(K)$ and $\mu
_{n}=\Omega c^{2}(K)/[h^{p}d(K)]$, where $\Omega$ is the
Lebesgue's measure of the support of $X_{t}$.
\end{theorem}

Theorem~\ref{th2} extends the results of Fan, Zhang and Zhang (\citeyear{FanZhaZha01}). We allow $X_{t}$
to be a vector and allow time series data. We do not assume that the
error $%
\varepsilon_{t}$ is independent of $X_{t}$ or the past history of $%
\{X_{t},Y_{t}\}$ so conditional heteroskedasticity in a time series context
is allowed. This is consistent with the empirical stylized fact of
volatility clustering for high frequency financial time series. We note that
the proof of the asymptotic normality of the GLR test in a time series
context is much more involved than in an i.i.d. context. It is
interesting to observe that the Wilks phenomena do not hold under
conditional heteroskedasticity because the factors $r(K)$ and $\mu_{n}$
involve the nuisance function $\sigma^{2}(X_{t})=\operatorname{var}(\varepsilon
_{t}|X_{t})$, which is unknown under $\mathbb{H}_{0}$. Conditional
homoskedasticity is required to ensure the Wilks phenomena. In this
case, $%
r(K)$ and $\mu_{n}$ are free of nuisance functions.

Like the $q_{n}$ test, we also consider the case of $p<4$. A higher
dimension $p$ could be allowed by suitably modifying factors $r(K)$ and
$\mu
_{n}$, which would depend on the unknown density $f(x)$ of $X_{t}$ and thus
are not free of nuisance functions, even under conditional
homoskedasticity.%

\section{Relative efficiency}\label{sec5}
We now compare the relative efficiency between the loss function test
$q_{n}$ and the GLR test $\lambda_{n}$ under
the class of local alternatives%
%
\begin{equation}
\mathbb{H}_{n}(a_{n})\dvtx g_{0}(X_{t})=g(X_{t},
\theta_{0})+a_{n}\delta(X_{t}), \label{eq5.1}
\end{equation}
where $\delta\dvtx\mathbb{R}\rightarrow\mathbb{R}$ is an unknown continuous
function with $E[\delta^{4}(X_{t})]\leq C$. The term $a_{n}\delta(X_{t})$
characterizes the departure of the model $g(X_{t},\theta_{0})$ from the
true function $g_{0}(X_{t})$ and the rate $a_{n}$ is the speed at which the
departure vanishes to $0$ as the sample size $n\rightarrow\infty$. For
notational simplicity, we have suppressed the dependence of $g_{0}(X_{t})$
on $n$ here. Without loss of generality, we assume that $\delta
(X_{t})$ is
uncorrelated with $X_{t}$, namely $E[\delta(X_{t})X_{t}]=0$.

\begin{theorem}[(Local power)]\label{th3}
Suppose Assumptions
\ref{assa1}--\ref{assa5} hold, $h\propto n^{-\omega}$ for $\omega\in(0,1/2p)$,
and $p<4$. Then \textup{(i)} under $\mathbb{H}_{n}(a_{n})$
with $a_{n}=n^{-1/2}h^{-p/4}$, we have
\[
\frac{s(K)q_{n}-\nu_{n}}{\sqrt{2\nu_{n}}}\stackrel{d}
{\longrightarrow}%
N(\psi,1)\qquad
\mbox{as }n\rightarrow\infty,
\]
where $\psi=\sigma^{2}E[\delta^{2}(X_{t})]/\sqrt{2b(K)\int
\sigma^{4}(x)\,dx}$, and $s(K)$ and $\nu
_{n}$ are as in Theorem~\ref{th1}. Suppose in addition $\operatorname{var}(\varepsilon
_{t}|X_{t})=\sigma^{2}$ almost surely. Then $\psi=\break E[\delta
^{2}(X_{t})]/ \sqrt{2b(K)\Omega}$.

\textup{(ii)} under $\mathbb{H}_{n}(a_{n})$ with $%
a_{n}=n^{-1/2}h^{-p/4}$, we have
\[
\frac{r(K)\lambda_{n}-\mu_{n}}{\sqrt{2\mu_{n}}}\stackrel{d} {%
\longrightarrow}N(\xi,1)\qquad
\mbox{as } n\rightarrow\infty,
\]
where $\xi=\sigma^{2}E[\delta^{2}(X_{t})]/[2\sqrt{2d(K)\int
\sigma^{4}(x)\,dx}]$, and $r(K)$ and $\mu_{n}$ are
as in Theorem~\ref{th2}. Suppose in addition $\operatorname{var}(\varepsilon_{t}|X_{t})=\sigma
^{2}$
almost surely. Then $\xi=E[\delta^{2}(X_{t})]/ [2\sqrt{%
2d(K)\Omega}]$.
\end{theorem}

When $X_{t}$ is a scalar (i.e., $p\!=\!1)$ and $h\!=\!n^{-2/9}$, the factor $%
a_{n}\!=\!n^{-1/2}h^{-p/4}=n^{-4/9}$ achieves the optimal rate in the sense of
Ingster (\citeyear{Ing93N1,Ing93N2,Ing93N3}).
Following a similar reasoning to Fan, Zhang and Zhang (\citeyear{FanZhaZha01}), we can show that the $q_{n}$ test can also detect local
alternatives with the optimal rate $n^{-2k/(4k+p)}$ in the sense of Ingster
(\citeyear{Ing93N1,Ing93N2,Ing93N3}), for the
function space $\mathcal{F}_{k}=\{ \delta\in
L^{2}\dvtx\int\delta^{(k)}(x)^{2}\,dx\leq C\}$. For $p=1$ and $k=2$,
this is
achieved by setting $h=n^{-2/9}$.

It is interesting to note that the noncentrality parameter $\psi$ of
the $%
q_{n}$ test is independent of the curvature parameter $D=d^{\prime
\prime
}(0)/2$ of the loss function $d(\cdot)$. This implies that all loss
functions satisfying Assumption~\ref{assa5} are asymptotically equally efficient
under $\mathbb{H}_{n}(a_{n})$ in terms of Pitman's efficiency criterion
[\citet{Pit79}, Chapter 7], although their shapes may be different.

While the $q_{n}$ and $\lambda_{n}$ tests achieve the same optimal
rate of
convergence in the sense of Ingster (\citeyear{Ing93N1,Ing93N2,Ing93N3}), Theorem~\ref{th4} below
shows that under the same set of regularity conditions [including the same
bandwidth $h$ and the same kernel $K(\cdot)$ for both tests],
$q_{n}$ is
asymptotically more efficient than $\lambda_{n}$ under $\mathbb{H}%
_{n}(a_{n})$.

\begin{theorem}[(Relative efficiency)]\label{th4}
Suppose
Assumptions~\ref{assa1}--\ref{assa5} hold, $h\propto n^{-\omega}$ for $\omega\in(0,1/2p)$
and $p<4$. Then Pitman's relative
efficiency of the $q_{n}$ test over the GLR $\lambda_{n}$ test under $\mathbb{H}_{n}(a_{n})$ with
$a_{n}=n^{-1/2}h^{-p/4}$ is given by
\begin{equation}
\operatorname{ARE}(q_{n}\dvtx\lambda_{n})= \biggl\{
\frac{\int [ 2\mathbf{K} (
\mathbf{u} ) -\int\mathbf{K} ( \mathbf{u}+\mathbf{v} )
\mathbf{K} ( \mathbf{v} ) \,d\mathbf{v} ] ^{2}\,d\mathbf{u}}{\int%
[ \int\mathbf{K} ( \mathbf{u}+\mathbf{v} ) \mathbf{K} (
\mathbf{v} ) \,d\mathbf{v} ] ^{2}\,d\mathbf{u}} \biggr\} ^{{1}/{(2-p\omega)}}, \label{eq5.2}
\end{equation}
where $\mathbf{K}(\mathbf{u})=\prod_{i=1}^{p}K(u_{i}),\mathbf
{u}%
=(u_{1},\ldots,u_{p})^{\prime}$. The asymptotic relative
efficiency $\operatorname{ARE}(q_{n}\dvtx\lambda_{n})$ is larger than 1 for any kernel satisfying
Assumption~\ref{assa4} and the condition of $K(\cdot)\leq1$.
\end{theorem}

Theorem~\ref{th4} holds under both conditional heteroskedasticity and conditional
homoskedasticity. It suggests that although the GLR $\lambda_{n}$
test is
a natural extension of the classical parametric LR test and is a generally
applicable nonparametric inference procedure with many appealing features,
it does not have the optimal power property of the classical LR test. In
particular, the GLR test is always asymptotically less efficient than the
loss function test under $\mathbb{H}_{n}(a_{n})$ whenever they use the same
kernel $K(\cdot)$ and the same bandwidth $h$, including the optimal kernel
and~the optimal bandwidth (if any) for the GLR test. The relative efficiency
gain of the loss function test over the GLR test holds even if the GLR
test $%
\lambda_{n}$ uses the true likelihood function. This result is in sharp
contrast to the classical LR test in a parametric setup, which is
asymptotically most powerful according to the Neyman--Pearson lemma.

Insight into the relative efficiency between $q_{n}$ and $\lambda_{n}$ can
be obtained by a Taylor expansion of the $\lambda_{n}$ statistic,%
%
\begin{equation}
\lambda_{n}=\frac{1}{2}\frac{\operatorname{SSR}_{0}-\operatorname{SSR}_{1}}{(\operatorname{SSR}_{1}/n)}+\mbox{
Remainder}, \label{eq5.3}
\end{equation}
where the remainder term is an asymptotically negligible higher order term
under $\mathbb{H}_{n}(a_{n})$. This is equivalent to use of the loss
function%
%
\begin{equation}
D(g,g_{\theta})= \bigl[ Y_{t}-g(X_{t},\theta)
\bigr] ^{2}- \bigl[ Y_{t}-g(X_{t}) \bigr]
^{2}. \label{eq5.4}
\end{equation}
When $g(X_{t})$ is close to $g(X_{t},\theta_{0})$, the first-order
term in
a Taylor expansion of $D(g,g_{\theta})$ around $g_{\theta_{0}}$ does not
vanish to 0 under $\mathbb{H}_{0}$. More specifically, the asymptotic
distribution of $\lambda_{n}$ is determined by the dominant term, 
%
\begin{eqnarray}\label{eq5.5}
\frac{1}{2}[\operatorname{SSR}_{0}-\operatorname{SSR}_{1}]&=&\frac{1}{2}
\Biggl[ \sum_{t=1}^{n}\hat{
\varepsilon%
}_{t}^{2}-\sum
_{t=1}^{n}\bigl[\hat{\varepsilon}_{t}-
\hat{m}_{h}(X_{t})\bigr]^{2} \Biggr]
\nonumber
\\[-8pt]
\\[-8pt]
\nonumber
&=&\sum
_{t=1}^{n}\hat{\varepsilon}_{t}
\hat{m}_{h}(X_{t})-\frac{1}{2}%
\sum
_{t=1}^{n}\hat{m}_{h}^{2}(X_{t}).
\end{eqnarray}
The first term in (\ref{eq5.5}) corresponds to the first-order term of a Taylor
expansion of (\ref{eq5.4}). It is a second-order $V$-statistic [\citet
{Ser80}], and
after demeaning, it can be approximated as a second-order degenerate $U$-statistic.
The second term in (\ref{eq5.5}) corresponds to the second-order
term of
a Taylor expansion of (\ref{eq5.4}). It is a third-order $V$-statistic and can be
approximated by a second-order degenerate $U$-statistic after demeaning.
These two degenerate $U$-statistics are of the same order of magnitude and
jointly determine the asymptotic distribution of $\lambda_{n}$. In
particular, the asymptotic variance of $\lambda_{n}$ is determined by the
variances of these two $U$-statistics and their covariance. In contrast,
under Assumption~\ref{assa5}, a Taylor expansion suggests that the asymptotic
distribution of the $q_{n}$ statistic is determined by
%
\begin{equation}
\hat{Q}_{n}=D\sum_{t=1}^{n}
\hat{m}_{h}^{2}(X_{t})+\mbox{Remainder}, \label{eq5.6}
\end{equation}
which corresponds to the second term in the expansion of $\operatorname{SSR}_{0}-\operatorname{SSR}_{1}$
in (\ref{eq5.5}). As it turns out, the asymptotic variance of this term alone is
always smaller than the variance of the difference of the two terms in
(\ref{eq5.5}). This leads to a more efficient test than the GLR test, as is
shown in
Theorem~\ref{th4}. We note that the first term in (\ref{eq5.5}), which causes an efficiency
loss for the GLR test relative to the $q_{n}$ test, is always present no
matter whether we use the bias-reduction device (i.e., estimating the
conditional mean of the estimated model residuals).

To assess the magnitude of the relative efficiency gain of the $q_{n}$ test
over the $\lambda_{n}$ test, we consider a few commonly used multiweight
kernels: the uniform, Epanechnikov, biweight and triweight kernels; see
Table~\ref{tab1} below. Suppose the bandwidth rate parameter $\omega=1/5$,
$2/9$,
respectively, in the univariate case (i.e., $p=1)$. The rate of $\omega=1/5$
gives the optimal bandwidth rate for estimating $g_{0}(\cdot)$, and the
rate of $\omega=2/9$ achieves the optimal convergence rate in the
sense of
Ingster (\citeyear{Ing93N1}, \citeyear{Ing93N2}, \citeyear{Ing93N3}).
Table~\ref{tab1} reports Pitman's asymptotic
relative efficiencies (ARE). The efficiency gain of using the $q_{n}$ test
is substantial, no matter if the bandwidth $h$ is of the order of
$n^{-1/5}$ or $%
n^{-2/9}$. Furthermore, there is little difference in the asymptotic
relative efficiency between the two choices of $h$. These are confirmed in
our simulation study below.

\begin{table}
\caption{Asymptotic relative efficiency of the loss function test
over the GLR test}\label{tab1}
\begin{tabular*}{\textwidth}{@{\extracolsep{\fill}}lcccc@{}}
\hline
& \textbf{Uniform} & \textbf{Epanechnikov} &
\textbf{Biweight} &
\textbf{Triweight}
\\
\hline
$K(u)$ &$\frac{1}{2}1(|u|\leq1)$ & $
\frac{3}{4} \bigl[ 1-u^{2} \bigr] 1(|u|\leq1)$ & $
\frac{15}{16} \bigl[ 1-u^{2} \bigr] ^{2}1(|u|\leq1)$ & $
\frac{35}{32} \bigl[ 1-u^{2} \bigr] ^{3}1(|u|\leq1)$
\\
$\operatorname{ARE}_{1}$ & $2.80$ & $2.04$ & $1.99$ & $1.98$
\\
$\operatorname{ARE}_{2}$ & $2.84$ & $2.06$ & $2.01$ & $1.99$
\\
\hline
\end{tabular*}
\tabnotetext[]{}{\textit{Note}: $\operatorname{ARE}$ denotes Pitman's asymptotic relative
efficiency of the loss function ${q}_{n}$ test to the
GLR
$\lambda_{n}$ test. $\operatorname{ARE}_{1}$ is for ${h=cn}^{-{ 1/5}}$ and ${\operatorname{ARE}}_{2}$
is for $h={cn}^{{ -2/9}}$, for ${0<c<\infty}$.}
\end{table}

We emphasize that Theorem~\ref{th4} does not imply that the GLR test should be
abandoned. Indeed, it is a natural extension of the classical LR test and
has many appealing features. It will remain as a useful, general
nonparametric inference procedure in practice.

While the relative efficiency of the loss function $q_{n}$ test over
the GLR
$\lambda_{n}$ test holds whenever the same bandwdith $h$ and the same
kernel $K(\cdot)$ are used, the choice of an optimal bandwidth remains an
important issue for each test. Theorems~\ref{th1}--\ref{th4} allow for a wide range of the
choices of $h$, but they do not provide a practical guidance on how to
choose $h$. In practice, a simple rule of thumb is to choose $%
h=S_{X}n^{-1/5} $ or $h=S_{X}n^{-2/9}$, where $S_{X}^{2}$ is the sample
variance of $\{X_{t}\}_{t=1}^{n}$. One could also choose a data-driven
bandwidth using a cross-validation procedure, that is, choose $h=\arg
\min_{c_{1}n^{-1/(p+4)}\leq h\leq c_{1}n^{-1/(p+4)}}\sum_{t=1}^{n} [
\hat{\varepsilon}_{t}-\hat{m}_{h,t}(X_{t}) ] ^{2}$ for some
prespecified\break \mbox{constants} $0<c_{1}<c_{2}<\infty$, where for each given
$t$, $\hat{m}_{h,t}(X_{t})$ is the leave-one-out estimator that is based on the
sample $\{ \hat{\varepsilon}_{s},X_{s}\}_{s=1,s\neq t}^{n}$. The bandwidth
based on cross-validation is asymptotically optimal for estimation in terms
of mean squared errors, but it may not be optimal for the $q_{n}$ and $%
\lambda_{n}$ tests. For testing problems, the central concern is the
Type I
error or Type II error, or both. Based on the Edgeworth expansion of the
asymptotic distribution of a test statistic, \citet{GaoGij08} show
that the choice of $h$ affects both Type I and Type II errors of a
closely related nonparametric test, and usually there exists a tradeoff
between Type I and Type II errors when choosing $h$. A sensible optimal rule
is to choose $h$ to maximize the power of a test given a significance level.
\citet{GaoGij08} derive the leading terms of the size and power
functions of their test statistic, and then choose a bandwidth to maximize
the power under a class of local alternatives similar to (\ref{eq5.1}) under a
controlled significance level, that is, to choose $h=\max_{h\in
B_{n}(\alpha
)}\beta_{n}(h)$, where $B_{n}(\alpha)=\{h\dvtx\alpha-c_{\min
}<\alpha
_{n}(h)<\alpha+c_{\min}\}$ for some prespecified small constant
$c_{\min
}\in(0,\alpha)$, and $\alpha_{n}(h)$ and $\beta_{n}(h)$ are the
size and
power functions of the nonparametric test. They then propose a data-driven
bandwidth in combination with a bootstrap and show that it works well in
finite samples. Unfortunately, Gao and Gijbels's (\citeyear{GaoGij08})
results cannot be
directly applied to either the $q_{n}$ or $\lambda_{n}$ test, because the
higher order terms of $\alpha_{n}(h)$ and $\beta_{n}(h)$ depend on the
form of test statistic, the DGP, the kernel $K$ and the bandwidth $h$, among
many other things. However, it is possible to extend their approach to
the $%
q_{n}$ and $\lambda_{n}$ tests to obain their optimal banwidths,
respectively. As the associated technicality is quite involved, we leave
this important problem for subsequent work. We note that \citet
{SunPhiJin08}, in a different context, also consider a data-driven
bandwidth by
minimizing a weighted average of the Type I and Type~II errors of a test,
namely choose $h=\arg\min_{h} ( \frac{w_{n}}{1+w_{n}}e_{n}^{I}+\frac
{1}{%
1+w_{n}}e_{n}^{\mathit{II}} )$, where $e_{n}^{I}$ and $e_{n}^{\mathit{II}}$ are the Type
I and Type II errors, respectively, and $w_{n}$ is a weight function that
reflects the relative importance of $e_{n}^{I}$ and $e_{n}^{\mathit{II}}$.

\section{Monte Carlo evidence}\label{sec6}
We now compare the finite sample
performance of the loss function $q_{n}$ test and the GLR $\lambda_{n}$
test. To examine the sizes of the tests, we consider the following null
linear regression model in a time series context:

\setcounter{dgp}{-1}
\begin{dgp}[(Linear regression)]\label{dgp0}
\[
\cases{
Y_{t}=1+X_{t}+\varepsilon
_{t},
\vspace*{2pt}\cr
X_{t}=0.5X_{t-1}+v_{t},
\vspace*{2pt}\cr
v_{t}\sim \mathrm{i.i.d.}\ N(0,1).}
\]
Here, $X_{t}$ is truncated within its two standard deviations. To examine
robustness of the tests, we consider a variety of distributions for the
error $\varepsilon_{t}$: (i)~$\varepsilon_{t}\sim \mathrm{i.i.d.}\ N(0,1)$,
(ii)~$%
\varepsilon_{t}\sim \mathrm{i.i.d.}$ Student-$t_{5}$, (iii) $\varepsilon
_{t}\sim
\mathrm{i.i.d.}$ $U[0,1]$, (iv)~$\varepsilon_{t}\sim \mathrm{i.i.d.}\ \ln N(0,1)$ and (v) $
\varepsilon_{t}\sim \mathrm{i.i.d.}\ \chi_{1}^{2}$, where the $\varepsilon_{t}$
in (iii)--(v) have been scaled to have mean $0$ and variance 1.
\end{dgp}

Because the asymptotic normal approximation for the $q_{n}$ and $\lambda
_{n} $ tests might not perform well in finite samples, we also use a
conditional bootstrap procedure based on the Wilks phenomena:

\textit{Step} 1: Obtain the parameter estimator $\hat{\theta}_{0}$
(e.g., OLS) of the null linear regression model, and the nonparametric
estimator $\hat{g}(X_{t})$.

\textit{Step} 2: Compute the $q_{n}$ statistic and the residual $%
\hat{\varepsilon}_{t}=Y_{t}-\hat{g}(X_{t})$ from the nonparametric model.

\textit{Step} 3: Conditionally on each $X_{t}$, draw a bootstrap
error $\varepsilon_{t}^{\ast}$ from the centered empirical
distribution of
$\hat{\varepsilon}_{t}$ and compute $Y_{t}^{\ast}=X_{t}^{\prime}\hat
{\theta%
}_{0}^{\ast}+\hat{\varepsilon}_{t}^{\ast}$. This forms a conditional
bootstrap sample $\{X_{t},Y_{t}^{\ast}\}_{t=1}^{n}$.

\textit{Step} 4: Use the conditional bootstrap sample $%
\{X_{t},Y_{t}^{\ast}\}_{t=1}^{n}$ to compute a bootstrap statistic $%
q_{n}^{\ast}$, using the same kernel $K(\cdot)$ and the same
bandwidth $h$
as in step 2.

\textit{Step} 5: Repeat steps 3 and 4 for a total of $B$ times,
where $B$ is a large number. We then obtain a collection of bootstrap test
statistics, $\{q_{nl}^{\ast}\}_{l=1}^{B}$.

\textit{Step} 6: Compute the bootstrap $P$ value $P^{\ast
}=B^{-1}\sum_{l=1}^{B}\mathbf{1}(q_{n}<q_{nl}^{\ast})$. Reject $\mathbb
{H}%
_{0}$ at a prespecified significance level $\alpha$ if and only if
$P^{\ast
}<\alpha$.

When conditional heteroskedasticity exists, we can modify step 2 by
using a
wild bootstrap for $\{ \hat{\varepsilon}_{t}^{\ast}\}$. If $X_{t}$ contains
lagged dependent variables, we can use a recursive simulation method;
see, for example, Franke, Kreiss and Mammen (\citeyear{FraKreMam02}). For space, we do not justify the
validity of the bootstrap here. Fan and Jiang [(\citeyear{FanJia07}),
Theorem 7] show the
consistency of the bootstrap for the GLR test in an i.i.d. context.
We could establish the consistency of the bootstrap for our loss function
test by following the approaches of \citet{FanJia07} and
\citet{GaoGij08}.

We consider two versions of the loss function test, one is to
standardize $%
\hat{Q}_{n}$ by $\hat{\sigma}_{n}^{2}=n^{-1}\operatorname{SSR}_{1}$, where $\operatorname{SSR}_{1}$
is the
sum of squared residuals of the nonparametric regression estimates.
This is
denoted as $q_{n}$. The other version is to standardize $\hat{Q}_{n}$
by $%
\hat{\sigma}_{n,0}^{2}=n^{-1}\operatorname{SSR}_{0}$, where $\operatorname{SSR}_{0}$ is the sum of squared
residuals of the null linear model. This is denoted as $q_{n}^{0}$. It is
expected that $q_{n}$ may be more powerful than $q_{n}^{0}$ in finite
samples under $\mathbb{H}_{A}$, because $\operatorname{SSR}_{0}$ is expected to be
significantly larger than $\operatorname{SSR}_{1}$ under $\mathbb{H}_{A}$. To
construct the
$q_{n}$ and $q_{n}^{0}$ tests, we choose the family of linex loss functions
in (\ref{eq3.12}), with $(\alpha,\beta)=(0,1),(0.2,1),(0.5,1)$ and $(1,1)$,
respectively; see Figure~\ref{fig1} for their shapes. The choice of $(\alpha,\beta)=(0,1)$ corresponds to the symmetric quadratic loss function, while
the degree of asymmetry of the loss function increases as $\alpha$
increases (the choice of $\beta$ has no impact on the $q_{n}$ tests).
Various choices of $(\alpha,\beta)$ thus allow us to examine sensitivity
of the power of the $q_{n}$ tests to the choices of the loss function.
Rather conveniently, when using the bootstrap procedure, there is no
need to
compute the centering and scaling factors for the $q_{n}$ and $q_{n}^{0}$
tests; it suffices to compare the statistic $q_{n}$ or $q_{n}^{0}$ with
their bootstrap counterparts. We choose $B=99$. The same bootstrap is used
for the GLR test $\lambda_{n}$.

To examine the power of the tests, we consider three nonlinear DGP's:

\begin{dgp}[(Quadratic regression)]\label{dgp1}
\[
Y_{t}=1+X_{t}+\theta X_{t}^{2}+
\varepsilon_{t}.
\]
\end{dgp}

\begin{dgp}[(Threshold regression)]\label{dgp2}
\[
Y_{t}=1+X_{t}\mathbf{1}(X_{t}>0)+(1+\theta
)X_{t}\mathbf{1}(X_{t}\leq0)+\varepsilon
_{t}.
\]
\end{dgp}

\begin{dgp}[(Smooth transition regression)]\label{dgp3}
\[
Y_{t}=1+X_{t}+\bigl[1-\theta F(X_{t})
\bigr]X_{t}+\varepsilon_{t},
\]
where $F(X_{t})=[1+\exp(-X_{t})]^{-1}$.
\end{dgp}

We consider various values for $\theta$ in each DGP to examine how the
power of the tests changes as the value of $\theta$ changes.

To examine sensitivity of all tests to the choices of $h$, we consider $
h=S_{X}n^{-\omega}$ for $\omega=\frac{2}{9}$ and $\frac{1}{5}$,
respectively, where $S_{X}$ is the sample standard deviation of $%
\{X_{t}\}_{t=1}^{n}$. These correspond to the optimal rate of
convergence in
the sense of Ingster (\citeyear{Ing93N1,Ing93N2,Ing93N3}) and the optimal rate of
estimation in terms of mean squared errors, respectively. The results are
similar. Here, we focus our discussion on the results with $h=S_{X}n^{-2/9}$,
as reported in Tables~\ref{tab2}--\ref{tab6}. The results with $h=S_{X}n^{-1/5}$ are
reported in Tables S.1--S.5 of the supplementary material. We use the uniform
kernel $K(z)=\frac{1}{2}\mathbf{1}(|z|\leq1)$ for all tests. We have also
used the biweight kernel, and the results are very similar (so, not reported
here).

\begin{sidewaystable}
\tablewidth=\textwidth
\caption{Empirical sizes of tests using asymptotic critical values}\label{tab2}
{\fontsize{7.3}{9.3}{\selectfont{
\begin{tabular*}{\textwidth}{@{\extracolsep{\fill}}ld{2.1}d{2.1}d{1.1}d{1.1}d{1.1}d{1.1}d{1.1}d{1.1}d{1.1}d{1.1}d{1.1}d{1.1}d{1.1}d{1.1}d{1.1}d{1.1}d{1.1}
d{1.1}@{}}
\hline
& \multicolumn{6}{c}{$\bolds{n=100}$} & \multicolumn{6}{c}{$\bolds{n=250}$}
& \multicolumn{6}{c@{}}{$\bolds{n=500}$} \\[-4pt]
& \multicolumn{6}{c}{\hrulefill} & \multicolumn{6}{c}{\hrulefill}
& \multicolumn{6}{c@{}}{\hrulefill} \\
& \multicolumn{2}{c}{$\bolds{{q}_{{n}}}$} &
\multicolumn{2}{c}{$\bolds{{q}_{{  n}}^{{ 0}}}$} & \multicolumn{2}{c}{${\operatorname{\mathbf{GLR}}}$} & \multicolumn{2}{c}{$\bolds{{q}_{{  n}}}$} &
\multicolumn{2}{c}{$%
\bolds{{  q}_{{  n}}^{{ 0}}}$} & \multicolumn{2}{c}{${\operatorname{\mathbf{GLR}}}$} & \multicolumn{2}{c}{$\bolds{{  q}_{{  n}}}$} &
\multicolumn{2}{c}{$\bolds{{  q}_{{  n}}^{{ 0}}}$} & \multicolumn{2}{c}{${\operatorname{\mathbf{GLR}}}$}
\\[-4pt]
& \multicolumn{2}{c}{\hrulefill}&\multicolumn{2}{c}{\hrulefill}&
\multicolumn{2}{c}{\hrulefill}&
\multicolumn{2}{c}{\hrulefill}&
\multicolumn{2}{c}{\hrulefill}&
\multicolumn{2}{c}{\hrulefill}&
\multicolumn{2}{c}{\hrulefill}&
\multicolumn{2}{c}{\hrulefill}&
\multicolumn{2}{c@{}}{\hrulefill}\\
\multicolumn{1}{@{}l}{$\bolds{{ (a,\beta)}}$} & \multicolumn{1}{c}{\textbf{{10\%}}} & \multicolumn{1}{c}{\textbf{{5\%}}} & \multicolumn{1}{c}{\textbf{{10\%}}} &
\multicolumn{1}{c}{\textbf{{5\%}}} & \multicolumn{1}{c}{\textbf{{10\%}}} & \multicolumn{1}{c}{\textbf{{5\%}}} & \multicolumn{1}{c}{\textbf{{10\%}}} &
\multicolumn{1}{c}{\textbf{{5\%}}} &
\multicolumn{1}{c}{\textbf{{10\%}}} & \multicolumn{1}{c}{\textbf{{5\%}}} & \multicolumn{1}{c}{\textbf{{10\%}}} & \multicolumn{1}{c}{\textbf{{5\%}}} &
\multicolumn{1}{c}{\textbf{{10\%}}}
& \multicolumn{1}{c}{\textbf{{5\%}}} & \multicolumn{1}{c}{\textbf{{10\%}}} & \multicolumn{1}{c}{\textbf{{5\%}}} & \multicolumn{1}{c}{\textbf{{10\%}}} &
\multicolumn{1}{c@{}}{\textbf{{5\%}}}
\\
\hline
&\multicolumn{18}{c}{DGP S.1: i.i.d. normal errors} \\
$(0.0,1.0)$ & 7.7 & 5.0 & 4.5 & 2.9 & 7.4 & 5.1 & 6.0 & 3.6 &  4.6 &  2.6 &  7.3 &  4.5 &
6.8 &  4.8 &  6.3 &  4.0 &  7.2 &
3.8 \\
 $(0.2,1.0)$ &  8.1 &  5.2 &  4.5 &
 3.0%
 &  7.4 &  5.1 &  6.1 &
3.8 &  4.7 &  2.5 &  7.3 &  4.5 &
6.6 &  4.9 &  6.2 &  4.1 &  7.2 &
3.8 \\
 $(0.5,1.0)$ &  8.4 &  5.2 &  4.9 &
 3.5%
 &  7.4 &  5.1 &  6.6 &
3.7 &  4.8 &  2.6 &  7.3 &  4.5 &
6.6 &  4.8 &  6.1 &  4.5 &  7.2 &
3.8 \\
 $(1.0,1.0)$ &  9.4 &  5.8 &  5.7 &
 4.3%
 &  7.4 &  5.1 &  6.8 &
4.3 &  5.2 &  3.0 &  7.3 &  4.5 &
6.9 &  5.0 &  6.0 &  4.6 &  7.2 &
3.8 \\
& \multicolumn{18}{c}{DGP S.2: i.i.d. Student-$t_{5}$ errors} \\
 $(0.0,1.0)$ &  6.8 &  4.3 &  4.0 &
 2.5%
 &  6.4 &  3.4 &  6.0 &  3.3 &
 4.4
&  2.4 &  4.9 &  2.7 &  5.4 &
3.1 &
 4.3 &  2.4 &  7.7 &  4.5 \\
 $(0.2,1.0)$ &  6.7 &  4.3 &  4.1 &
 2.6%
 &  6.4 &  3.4 &  5.8 &  3.8 &
 4.4
&  2.4 &  4.9 &  2.7 &  5.4 &
3.2 &
 4.6 &  2.5 &  7.7 &  4.5 \\
 $(0.5,1.0)$ &  6.9 &  4.5 &  4.3 &
 2.6%
 &  6.4 &  3.4 &  5.9 &  3.9 &
 4.3
&  2.7 &  4.9 &  2.7 &  5.6 &
3.3 &
 4.8 &  2.6 &  7.7 &  4.5 \\
 $(1.0,1.0)$ &  8.5 &  5.1 &  5.0 &
 3.0%
 &  6.4 &  3.4 &  6.4 &  4.4 &
 5.0
&  3.4 &  4.9 &  2.7 &  6.2 &
3.5 &
 5.1 &  3.0 &  7.7 &  4.5 \\
& \multicolumn{18}{c}{DGP S.3: i.i.d. uniform errors} \\
 $(0.0,1.0)$ &  7.1 &  5.2 &  4.3 &
 2.7%
 &  6.2 &  3.5 &  6.8 &  4.2 &
 5.5
&  2.9 &  6.6 &  4.0 &  6.4 &
4.2 &
 5.6 &  3.5 &  6.3 &  3.2 \\
 $(0.2,1.0)$ &  7.1 &  5.4 &  4.5 &
 2.5%
 &  6.2 &  3.5 &  6.8 &  4.4 &
 5.4
&  2.9 &  6.6 &  4.0 &  6.2 &
4.2 &
 5.3 &  3.5 &  6.3 &  3.2 \\
 $(0.5,1.0)$ &  7.4 &  5.4 &  4.9 &
 2.7%
 &  6.2 &  3.5 &  7.0 &  4.5 &
 5.6
&  3.1 &  6.6 &  4.0 &  6.2 &
4.5 &
 5.5 &  3.6 &  6.3 &  3.2 \\
 $(1.0,1.0)$ &  9.0 &  6.0 &  5.7 &
 3.7%
 &  6.2 &  3.5 &  7.1 &  5.1 &
 5.9
&  3.4 &  6.6 &  4.0 &  6.7 &
4.6 &
 6.1 &  3.7 &  6.3 &  3.2 \\
&\multicolumn{18}{c}{DGP S.4: i.i.d. log-normal
errors}
\\
 $(0.0,1.0)$ &  9.4 &  7.2 &  6.8 &
 4.3%
 &  8.5 &  6.6 &  6.1 &  4.2 &
 4.8
&  2.8 &  6.7 &  3.6 &  6.9 &
5.1 &
 6.3 &  4.4 &  6.9 &  4.1 \\
 $(0.2,1.0)$ &  9.9 &  7.9 &  7.6 &
 5.2%
 &  8.5 &  6.6 &  6.5 &  4.7 &
 4.9
&  3.3 &  6.7 &  3.6 &  7.5 &
5.6 &
 6.6 &  4.5 &  6.9 &  4.1 \\
 $(0.5,1.0)$ &  10.4 &  8.7 &  8.2 &
6.5 &  8.5 &  6.6 &  8.1 &  4.8 &
5.9 &  4.3 &  6.7 &  3.6 &  7.9 &
5.9 &  7.1 &  5.3 &  6.9 &  4.1 \\
 $(1.0,1.0)$ &  12.4 &  10.2 &  9.5 &
7.9 &  8.5 &  6.6 &  9.2 &  6.9 &
7.6 &  5.5 &  6.7 &  3.6 &  9.5 &
7.1 &  8.3 &  6.2 &  6.9 &  4.1 \\
&\multicolumn{18}{c}{DGP S.5: i.i.d. chi-square errors} \\
 $(0.0,1.0)$ &  7.6 &  5.9 &  5.6 &
 3.5%
 &  7.3 &  5.2 &  6.3 &  4.0 &
 4.6
&  2.8 &  6.3 &  3.1 &  5.2 &
3.5 &
 4.4 &  2.8 &  5.4 &  2.9 \\
 $(0.2,1.0)$ &  8.1 &  6.3 &  6.0 &
 3.8%
 &  7.3 &  5.2 &  6.5 &  4.2 &
 5.1
&  3.0 &  6.3 &  3.1 &  5.2 &
3.7 &
 4.9 &  3.1 &  5.4 &  2.9 \\
 $(0.5,1.0)$ &  8.8 &  7.2 &  7.0 &
 4.6%
 &  7.3 &  5.2 &  7.0 &  4.8 &
 5.6
&  3.7 &  6.3 &  3.1 &  5.5 &
3.9 &
 5.3 &  3.3 &  5.4 &  2.9 \\
 $(1.0,1.0)$ &  10.8 &  9.2 &  8.9 &
6.2 &  7.3 &  5.2 &  7.8 &  5.6 &
6.3 &  4.8 &  6.3 &  3.1 &  6.3 &
4.6 &  5.9 &  4.1 &  5.4 &  2.9 \\
\hline
\end{tabular*}}}}
\tabnotetext[]{}{
\fontsize{7.2}{9.2}{\selectfont{\textit{Notes}:
(i) 1000 iterations;
(ii) $\operatorname{GLR}$, the generalized likelihood ratio test, ${q}_{n}$ and ${q}_{{n}}^{{0}}$, loss function-based tests;
(iii) ${q}_{n}$ is standardized by ${\operatorname{SSR}}_{1}$,\ the sum of squared residuals of the nonparametric
regression
estimates, and $q_{{n}}^{{ 0}}$ is
standardized by ${\operatorname{SSR}}_{0}$, the sum of squared
residuals
of the null linear model;
(iv) The uniform kernel is used for ${\operatorname{GLR}}$,
${q}_{n}$ and ${q}_{{  n}}^{{ 0}}$;
the bandwidth ${h=S}_{{X}}{  n}^{-2/9}$,
where ${  S}_{{  X}}$ is the sample standard
deviation
of $\{X_{t}\}_{t=1}^{n}$;
(v) The ${  q}_{{  n}}$ tests are based on
the
linex loss function: ${d(z)=\frac{\beta}{\alpha^{2}}[\exp(
\alpha
z)-1-\alpha z];}$
(vi) $ Y_{t}= 1+X_{t} +\varepsilon
_{t}$, ${X}_{t}={ 0.5X}_{t-1}+{ v}_{t}$,
${v}_{t}\sim \mathrm{i.i.d.}\ N(0,1)$, where DGP S.1: $\varepsilon_{i}\thicksim \mathrm{i.i.d.}\ N(0,1)$;
DGP S.2: ${\varepsilon}_{i}
\thicksim \mathrm{i.i.d.}$ Student-$t_{{ 5}}$;
DGP S.3: ${\varepsilon}_{i} \thicksim \mathrm{i.i.d.}\ U[0,1]$;
DGP S.4: $\varepsilon_{i}  \thicksim \mathrm{i.i.d.}\ \log N(0,1)$;
DGP S.5: ${ \varepsilon}_{i} \thicksim \mathrm{i.i.d.}\ \chi_{{
1}}^{2}$.}}}
\end{sidewaystable}

\begin{sidewaystable}
\tablewidth=\textwidth
\caption{Empirical sizes of tests using bootstrap critical values}\label{tab3}
{\fontsize{7.3}{9.3}{\selectfont{\begin{tabular*}{\textwidth}{@{\extracolsep{\fill}}ld{2.1}d{1.1}
d{2.1}d{1.1}
d{2.1}d{1.1}
d{2.1}d{1.1}
d{2.1}d{1.1}
d{2.1}d{1.1}
d{2.1}d{1.1}
d{2.1}d{1.1}
d{2.1}d{1.1}
d{2.1}d{1.1}@{}}
\hline
& \multicolumn{6}{c}{$\bolds{n=100}$} & \multicolumn{6}{c}{$\bolds{n=250}$}
& \multicolumn{6}{c@{}}{$\bolds{n=500}$} \\[-4pt]
& \multicolumn{6}{c}{\hrulefill} & \multicolumn{6}{c}{\hrulefill}
& \multicolumn{6}{c@{}}{\hrulefill} \\
& \multicolumn{2}{c}{$\bolds{{q}_{{n}}}$} &
\multicolumn{2}{c}{$\bolds{{q}_{{  n}}^{{ 0}}}$} & \multicolumn{2}{c}{${\operatorname{\mathbf{GLR}}}$} & \multicolumn{2}{c}{$\bolds{{q}_{{  n}}}$} &
\multicolumn{2}{c}{$%
\bolds{{  q}_{{  n}}^{{ 0}}}$} & \multicolumn{2}{c}{${\operatorname{\mathbf{GLR}}}$} & \multicolumn{2}{c}{$\bolds{{  q}_{{  n}}}$} &
\multicolumn{2}{c}{$\bolds{{  q}_{{  n}}^{{ 0}}}$} & \multicolumn{2}{c}{${\operatorname{\mathbf{GLR}}}$}
\\[-4pt]
& \multicolumn{2}{c}{\hrulefill}&\multicolumn{2}{c}{\hrulefill}&
\multicolumn{2}{c}{\hrulefill}&
\multicolumn{2}{c}{\hrulefill}&
\multicolumn{2}{c}{\hrulefill}&
\multicolumn{2}{c}{\hrulefill}&
\multicolumn{2}{c}{\hrulefill}&
\multicolumn{2}{c}{\hrulefill}&
\multicolumn{2}{c@{}}{\hrulefill}\\
\multicolumn{1}{@{}l}{$\bolds{{ (a,\beta)}}$} & \multicolumn{1}{c}{\textbf{{10\%}}} & \multicolumn{1}{c}{\textbf{{5\%}}} & \multicolumn{1}{c}{\textbf{{10\%}}} &
\multicolumn{1}{c}{\textbf{{5\%}}} & \multicolumn{1}{c}{\textbf{{10\%}}} & \multicolumn{1}{c}{\textbf{{5\%}}} & \multicolumn{1}{c}{\textbf{{10\%}}} &
\multicolumn{1}{c}{\textbf{{5\%}}} &
\multicolumn{1}{c}{\textbf{{10\%}}} & \multicolumn{1}{c}{\textbf{{5\%}}} & \multicolumn{1}{c}{\textbf{{10\%}}} & \multicolumn{1}{c}{\textbf{{5\%}}} &
\multicolumn{1}{c}{\textbf{{10\%}}}
& \multicolumn{1}{c}{\textbf{{5\%}}} & \multicolumn{1}{c}{\textbf{{10\%}}} & \multicolumn{1}{c}{\textbf{{5\%}}} & \multicolumn{1}{c}{\textbf{{10\%}}} &
\multicolumn{1}{c@{}}{\textbf{{5\%}}}
\\
\hline
&\multicolumn{18}{c}{DGP S.1: i.i.d. normal errors} \\
 $(0.0,1.0)$ &  10.0 &  5.4 &  10.0 &
5.3 &  10.4 &  4.4 &  11.1 &  6.1 &
11.8 &  5.9 &  11.5 &  5.0 &  11.4 &
 6.1 &  11.6 &  6.1 &  10.9 &
6.5
\\
 $(0.2,1.0)$ &  10.0 &  5.4 &  10.0 &
5.6 &  10.4 &  4.4 &  11.8 &  6.2 &
12.1 &  6.3 &  11.5 &  5.0 &  11.4 &
 6.1 &  11.4 &  6.0 &  10.9 &
6.5
\\
 $(0.5,1.0)$ &  9.4 &  4.9 &  9.7 &
 4.9%
 &  10.4 &  4.4 &  11.4 &  5.9 &
11.7 &  5.8 &  11.5 &  5.0 &  11.4 &
 6.4 &  11.6 &  6.2 &  10.9 &
6.5
\\
 $(1.0,1.0)$ &  9.2 &  4.6 &  9.7 &
 4.5%
 &  10.4 &  4.4 &  11.3 &  5.8 &
11.2 &  5.5 &  11.5 &  5.0 &  11.7 &
 6.4 &  12.2 &  6.1 &  10.9 &
6.5
\\
&\multicolumn{18}{c}{DGP S.2: i.i.d. Student-$t_{5}$ errors} \\
 $(0.0,1.0)$ &  9.1 &  4.1 &  8.8 &
 4.4%
 &  9.5 &  3.8 &  8.8 &  4.3 &
 8.9
&  4.5 &  9.3 &  3.9 &  10.1 &
4.7
&  10.7 &  5.2 &  11.5 &  5.7 \\
 $(0.2,1.0)$ &  8.8 &  4.2 &  8.7 &
 4.2%
 &  9.5 &  3.8 &  8.8 &  4.3 &
 9.0
&  4.6 &  9.3 &  3.9 &  10.6 &
5.0
&  10.8 &  4.8 &  11.5 &  5.7 \\
 $(0.5,1.0)$ &  9.4 &  4.0 &  8.7 &
 4.4%
 &  9.5 &  3.8 &  9.1 &  4.5 &
 8.8
&  4.5 &  9.3 &  3.9 &  10.5 &
5.1
&  10.3 &  5.0 &  11.5 &  5.7 \\
 $(1.0,1.0$) &  9.8 &  3.9 &  10.2 &
4.5 &  9.5 &  3.8 &  8.8 &  5.0 &
9.2 &  4.7 &  9.3 &  3.9 &  10.6 &
5.4 &  10.5 &  5.2 &  11.5 &  5.7 \\
&\multicolumn{18}{c}{DGP S.3: i.i.d. uniform errors} \\
 $(0.0,1.0)$ &  10.3 &  5.0 &  10.1 &
5.3 &  9.1 &  4.2 &  11.1 &  5.8 &
11.2 &  5.8 &  10.9 &  6.5 &  9.6 &
6.0 &  9.5 &  6.1 &  10.6 &  5.3 \\
 $(0.2,1.0)$ &  10.3 &  5.1 &  10.3 &
5.5 &  9.1 &  4.2 &  11.1 &  5.8 &
10.9 &  5.7 &  10.9 &  6.5 &  9.6 &
6.0 &  9.5 &  6.0 &  10.6 &  5.3 \\
 $(0.5,1.0)$ &  10.8 &  5.1 &  10.7 &
5.4 &  9.1 &  4.2 &  11.2 &  5.6 &
11.2 &  5.7 &  10.9 &  6.5 &  9.4 &
6.0 &  9.4 &  6.0 &  10.6 &  5.3 \\
 $(1.0,1.0)$ &  10.3 &  5.5 &  10.6 &
5.2 &  9.1 &  4.2 &  10.7 &  5.7 &
11.0 &  5.8 &  10.9 &  6.5 &  9.5 &
6.0 &  9.4 &  6.2 &  10.6 &  5.3 \\
&\multicolumn{18}{c}{DGP S.4: i.i.d. log-normal errors}
\\
 $(0.0,1.0)$ &  11.0 &  5.8 &  10.9 &
6.1 &  10.2 &  5.5 &  9.7 &  4.6 &
9.6 &  5.1 &  10.8 &  4.5 &  10.3 &
5.3 &  10.4 &  5.4 &  11.0 &  5.9 \\
 $(0.2,1.0)$ &  10.7 &  5.9 &  10.6 &
6.0 &  10.2 &  5.5 &  10.0 &  4.7 &
9.5 &  4.9 &  10.8 &  4.5 &  9.9 &
5.3 &  10.3 &  5.4 &  11.0 &  5.9 \\
 $(0.5,1.0)$ &  10.9 &  5.8 &  10.8 &
5.6 &  10.2 &  5.5 &  9.5 &  4.7 &
9.4 &  4.7 &  10.8 &  4.5 &  9.9 &
5.6 &  10.0 &  5.4 &  11.0 &  5.9 \\
 $(1.0,1.0)$ &  10.7 &  5.8 &  11.0 &
5.9 &  10.2 &  5.5 &  9.7 &  4.5 &
9.5 &  4.4 &  10.8 &  4.5 &  10.1 &
5.6 &  10.4 &  5.4 &  11.0 &  5.9 \\
&\multicolumn{18}{c}{DGP S.5: i.i.d. chi-square errors}\\
 $(0.0,1.0)$ &  9.4 &  4.4 &  9.2 &
 4.4%
 &  9.7 &  5.3 &  10.2 &  4.8 &
 10.2%
 &  4.9 &  9.3 &  4.6 &  7.9 &
 3.7
&  8.2 &  3.7 &  9.0 &  4.1 \\
 $(0.2,1.0)$ &  9.3 &  4.5 &  9.0 &
 4.3%
 &  9.7 &  5.3 &  10.2 &  4.9 &
 10.0%
 &  4.9 &  9.3 &  4.6 &  7.7 &
 3.7
&  8.0 &  3.7 &  9.0 &  4.1 \\
 $(0.5,1.0)$ &  9.2 &  4.9 &  9.1 &
 4.8%
 &  9.7 &  5.3 &  10.2 &  4.7 &
 10.2%
 &  4.9 &  9.3 &  4.6 &  7.8 &
 3.7
&  7.9 &  3.7 &  9.0 &  4.1 \\
 $(1.0,1.0)$ &  8.7 &  5.2 &  8.5 &
 4.8%
 &  9.7 &  5.3 &  10.3 &  3.7 &
 9.9
&  3.6 &  9.3 &  4.6 &  7.8 &
3.6 &
 7.7 &  3.7 &  9.0 &  4.1 \\
 \hline
\end{tabular*}}}}
\tabnotetext[]{}{
{\fontsize{7.2}{9.2}{\selectfont{
\textit{Notes}:
(i) 1000 iterations;
(ii) ${  \operatorname{GLR}}$, the generalized likelihood ratio
test, ${  q}_{n}$ and ${  q}_{{  n}}^{{
 0}}$, loss function-based tests;
(iii) ${  q}_{n}$ is standardized by ${\operatorname{SSR}}_{1}$, the sum of squared residuals of the nonparametric
regression
estimates, and ${q}_{{  n}}^{{ 0}}$ is
standardized by ${\operatorname{SSR}}_{0}$, the sum of squared
residuals
of the null linear model;
(iv) The uniform kernel is used for ${\operatorname{GLR}}$,
$q_{n}$ and ${  q}_{{  n}}^{{ 0}}$;
the bandwidth $  h=S_{{  X}}{  n}^{-2/9}$, where ${  S}_{{  X}}$ is the sample standard
deviation
of $\{{  X}_{t}\}_{t=1}^{n}$;
(v) The ${  q}_{{  n}}$ tests are based on
the
linex loss function: $d(z)=\frac{\beta}{\alpha^{2}}[\exp(
\alpha
z)-1-\alpha z];$
(vi) ${  Y}_{t} =1+X_{t} +\varepsilon
_{t}$, ${X}_{t}=0.5X_{t-1} +v_{t}$,
${v}_{t}\sim \mathrm{i.i.d.}\ N(0,1)$,
where DGP S.1: $\varepsilon_{i} \thicksim \mathrm{i.i.d.}\ N(0,1)$;
DGP S.2: ${\varepsilon}_{i}\thicksim \mathrm{i.i.d.}\  \mathrm{Student}$-$t_{{ 5}}$;
DGP S.3: ${\varepsilon}_{i} \thicksim \mathrm{i.i.d.}\  U[0,1]$;
DGP S.4: $\varepsilon_{i}  \thicksim \mathrm{i.i.d.}\ \log N(0,1);$
DGP S.5: $\varepsilon_{i} \thicksim \mathrm{i.i.d.}\ \chi_{{
1}}^{2}$.}}}}
\end{sidewaystable}

Tables~\ref{tab2} and~\ref{tab3} report the empirical rejection rates of the tests
under $%
\mathbb{H}_{0}$ (DGP~\ref{dgp0}) at the 10\% and 5\% levels, using both asymptotic
and bootstrap critical values, respectively. We first examine the size
of the
tests using asymptotic critical values, with $n=100,250$ and $500$,
respectively. Table~\ref{tab2} shows that all tests, $\lambda_{n},q_{n}$ and $
q_{n}^{0}$, have reasonable sizes in finite samples, and they are
robust to
various error distributions, but they all show some underrejection,
particularly at the 10\% level. The $q_{n}$ and $\lambda_{n}$ tests have
similar sizes in most cases, whereas $q_{n}^{0}$ shows a bit more
underrejection. Overall, the sizes of the $q_{n},q_{n}^{0}$ and $\lambda
_{n} $ tests display some underrejections in most cases in finite samples,
but they are not unreasonable.

\begin{sidewaystable}
\tabcolsep=0pt
\tablewidth=\textwidth
\caption{Empirical powers of tests using bootstrap critical values}\label{tab4}
{\fontsize{8}{10}{\selectfont{\begin{tabular*}{\textwidth}{@{\extracolsep{4in minus 4in}}lcd{3.1}d{3.1}
d{3.1}d{3.1}
d{3.1}d{3.1}
d{3.1}d{3.1}
d{3.1}d{3.1}
d{3.1}d{3.1}
d{3.1}d{3.1}
d{3.1}d{3.1}
d{3.1}d{3.1}@{}}
\hline
&& \multicolumn{6}{c}{$\bolds{n=100}$} & \multicolumn{6}{c}{$\bolds{n=250}$}
& \multicolumn{6}{c@{}}{$\bolds{n=500}$} \\[-4pt]
& &\multicolumn{6}{c}{\hrulefill} & \multicolumn{6}{c}{\hrulefill}
& \multicolumn{6}{c@{}}{\hrulefill} \\
& &\multicolumn{2}{c}{$\bolds{{q}_{{n}}}$} &
\multicolumn{2}{c}{$\bolds{{q}_{{  n}}^{{ 0}}}$} & \multicolumn{2}{c}{${\operatorname{\mathbf{GLR}}}$} & \multicolumn{2}{c}{$\bolds{{q}_{{  n}}}$} &
\multicolumn{2}{c}{$%
\bolds{{  q}_{{  n}}^{{ 0}}}$} & \multicolumn{2}{c}{${\operatorname{\mathbf{GLR}}}$} & \multicolumn{2}{c}{$\bolds{{  q}_{{  n}}}$} &
\multicolumn{2}{c}{$\bolds{{  q}_{{  n}}^{{ 0}}}$} & \multicolumn{2}{c@{}}{${\operatorname{\mathbf{GLR}}}$}
\\[-4pt]
& &\multicolumn{2}{c}{\hrulefill}&\multicolumn{2}{c}{\hrulefill}&
\multicolumn{2}{c}{\hrulefill}&
\multicolumn{2}{c}{\hrulefill}&
\multicolumn{2}{c}{\hrulefill}&
\multicolumn{2}{c}{\hrulefill}&
\multicolumn{2}{c}{\hrulefill}&
\multicolumn{2}{c}{\hrulefill}&
\multicolumn{2}{c@{}}{\hrulefill}\\
\multicolumn{1}{@{}l}{$\bolds{(a,\beta)}$} &\multicolumn{1}{c}{$\bolds{\theta}$}& \multicolumn{1}{c}{\textbf{{10\%}}} & \multicolumn{1}{c}{\textbf{{5\%}}} & \multicolumn{1}{c}{\textbf{{10\%}}} &
\multicolumn{1}{c}{\textbf{{5\%}}} & \multicolumn{1}{c}{\textbf{{10\%}}} & \multicolumn{1}{c}{\textbf{{5\%}}} & \multicolumn{1}{c}{\textbf{{10\%}}} &
\multicolumn{1}{c}{\textbf{{5\%}}} &
\multicolumn{1}{c}{\textbf{{10\%}}} & \multicolumn{1}{c}{\textbf{{5\%}}} & \multicolumn{1}{c}{\textbf{{10\%}}} & \multicolumn{1}{c}{\textbf{{5\%}}} &
\multicolumn{1}{c}{\textbf{{10\%}}}
& \multicolumn{1}{c}{\textbf{{5\%}}} & \multicolumn{1}{c}{\textbf{{10\%}}} & \multicolumn{1}{c}{\textbf{{5\%}}} & \multicolumn{1}{c}{\textbf{{10\%}}} &
\multicolumn{1}{c@{}}{\textbf{{5\%}}}
\\
\hline
&&\multicolumn{18}{c}{DGP P.1: Quadratic regression} \\
 $(0.0,1.0)$ &  0.1 &  22.6 &  12.8 &
23.0 &  12.4 &  20.2 &
12.4 &
 38.4 &  27.2 &  39.2 &  27.8 &
 29.4
&  20.4 &  62.0 &  50.6 &  62.2 &
50.8 &  47.4 &  33.6 \\
&  0.2 &  53.2 &  39.8 &  53.6 &
 39.6 &  45.2 &  34.8 &
90.6 &  83.4 &  91.0 &  83.4 &  80.4 &
 70.6 &  99.4 &  99.0 &  99.4 &
 99.0
&  97.4 &  95.6 \\
&  0.3 &  85.2 &  77.0 &  85.2 &
 77.2 &  76.8 &  65.8 &
99.6 &  99.2 &  99.6 &  99.2 &  98.8 &
 97.2 &  100.0 &  100.0 &  100.0 &
100.0 &  100.0 &  100.0 \\
&  0.5 &  99.4 &  98.4 &  99.4 &
 98.8%
 &  98.6 &  96.6 &  100.0 &  100.0 &
 100.0 &  100.0 &  100.0 &  100.0 &
100.0 &  100.0 &  100.0 &  100.0 &
100.0 &
 100.0 \\
&  1.0 &  100.0 &  100.0 &  100.0 &
100.0 &  100.0 &  100.0 &  100.0 &
100.0 &
 100.0 &  100.0 &  100.0 &  100.0 &
100.0 &  100.0 &  100.0 &  100.0 &
100.0 &
 100.0 \\ [3pt]
 $(0.2,1.0)$ &  0.1 &  23.4 &  13.6 &
23.6 &  13.0 &  20.2 &  12.4 &  39.8 &
 28.0 &  40.0 &  28.4 &  29.4 &
 20.4
&  63.8 &  51.8 &  64.0 &  52.0 &
47.4 &  33.6 \\
&  0.2 &  55.8 &  41.2 &  55.8 &
 41.0%
 &  45.2 &  34.8 &  91.0 &  84.0 &
91.4 &  83.8 &  80.4 &  70.6 &  99.4 &
 99.2 &  99.4 &  99.2 &  97.4 &
 95.6
\\
&  0.3 &  86.2 &  78.4 &  86. &
 78.6
&  76.8 &  65.8 &  99.8 &  99.2 &
99.6 &  99.2 &  98.8 &  97.2 &  100.0
&
 100.0 &  100.0 &  100.0 &  100.0 &
100.0 \\
&  0.5 &  99.4 &  98.6 &  99.4 &
 99.2%
 &  98.6 &  96.6 &  100.0 &  100.0 &
 100.0 &  100.0 &  100.0 &  100.0 &
100.0 &  100.0 &  100.0 &  100.0 &
100.0 &
 100.0 \\
&  1.0 &  100.0 &  100.0 &  100.0 &
100.0 &  100.0 &  100.0 &  100.0 &
100.0 &
 100.0 &  100.0 &  100.0 &  100.0 &
100.0 &  100.0 &  100.0 &  100.0 &
100.0 &
 100.0 \\[3pt]
 $(0.5,1.0)$ &  0.1 &  24.8 &  14.4 &
24.2 &  14.6 &  20.2 &  12.4 &  42.6 &
 30.0 &  41.4 &  29.2 &  29.4 &
 20.4
&  64.8 &  53.4 &  64.6 &  53.6 &
47.4 &  33.6 \\
&  0.2 &  58.0 &  43.8 &  59.0 &
 43.4%
 &  45.2 &  34.8 &  91.4 &  85.2 &
92.0 &  85.8 &  80.4 &  70.6 &  99.6 &
 99.4 &  99.6 &  99.4 &  97.4 &
 95.6
\\
&  0.3 &  87.0 &  80.0 &  86.6 &
 80.2%
 &  76.8 &  65.8 &  99.8 &  99.4 &
99.6 &  99.4 &  98.8 &  97.2 &  100.0
&
 100.0 &  100.0 &  100.0 &  100.0 &
100.0 \\
&  0.5 &  99.4 &  98.8 &  99.4 &
 99.2%
 &  98.6 &  96.6 &  100.0 &  100.0 &
 100.0 &  100.0 &  100.0 &  100.0 &
100.0 &  100.0 &  100.0 &  100.0 &
100.0 &
 100.0 \\
&  1.0 &  100.0 &  100.0 &  100.0 &
100.0 &  100.0 &  100.0 &  100.0 &
100.0 &
 100.0 &  100.0 &  100.0 &  100.0 &
100.0 &  100.0 &  100.0 &  100.0 &
100.0 &
 100.0 \\[3pt]
 $(1.0,1.0)$ &  0.1 &  27.6 &  16.2 &
26.2 &  18.0 &  20.2 &  12.4 &  44.6 &
 32.4 &  45.2 &  31.4 &  29.4 &
 20.4
&  66.4 &  66.0 &  66.0 &  56.6 &
47.4 &  33.6 \\
&  0.2 &  60.8 &  46.4 &  62.2 &
 46.4%
 &  45.2 &  34.8 &  92.2 &  87.0 &
92.6 &  88.2 &  80.4 &  70.6 &  99.6 &
 99.4 &  99.8 &  99.4 &  97.4 &
 95.6
\\
&  0.3 &  88.8 &  81.6 &  88.4 &
 80.6%
 &  76.8 &  65.8 &  99.8 &  99.6 &
99.6 &  97.2 &  98.8 &  97.2 &  100.0
&
 100.0 &  100.0 &  100.0 &  100.0 &
100.0 \\
&  0.5 &  99.6 &  99.0 &  99.6 &
 99.2%
 &  98.6 &  96.6 &  100.0 &  100.0 &
 100.0 &  100.0 &  100.0 &  100.0 &
100.0 &  100.0 &  100.0 &  100.0 &
100.0 &
 100.0 \\
&  1.0 &  100.0 &  100.0 &  100.0 &
100.0 &  100.0 &  100.0 &  100.0 &
100.0 &
 100.0 &  100.0 &  100.0 &  100.0 &
100.0 &  100.0 &  100.0 &  100.0 &
100.0 &
 100.0 \\
 \hline
\end{tabular*}}}}
\tabnotetext[]{}{\textit{Notes}:
(i) 1000 iterations;
(ii) $\operatorname{GLR}$, the generalized likelihood ratio
test, ${  q}_{{  n}}$ and ${  q}_{{
n}}^{%
{ 0}}$,  loss function-based tests;
(iii) ${  q}_{{  n}}$ is standardized by
${  \operatorname{SSR}}_{1}$, the sum of squared residuals of the
nonparametric regression estimates, and ${  q}_{{  n}}^{{
 0}
}$ is standardized by ${  \operatorname{SSR}}_{0}$, the sum of
squared residuals of the null linear model;
(iv)~The uniform kernel is used for ${  \operatorname{GLR}}$, ${q}_{{  n}}$ and ${  q}_{{  n}}^{{ 0}}$;
the bandwidth ${  h=S}_{{  X}}{  n}^{-2/9}$, where ${  S}_{{  X}}$ is the sample standard
deviation
of $\{{  X}_{t}\}_{t=1}^{n}$;
(v) The ${  q}_{{  n}}$ tests are based on
the
linex loss function: $d(z)=\frac{\beta}{\alpha^{2}}[\exp(
\alpha
z)-1-\alpha z];$
(vi) DGP P.1, ${  Y}_{t} =1+X_{t}+\theta
X%
_{t}^{2} +\varepsilon_{t}$,  where $\{
\varepsilon_{i}\} \sim \mathrm{i.i.d.}\ N(0,1)$.}
\end{sidewaystable}

\begin{sidewaystable}
\tabcolsep=0pt
\tablewidth=\textwidth
\caption{Empirical powers of tests using bootstrap critical values}\label{tab5}
{\fontsize{7.1}{9.1}{\selectfont{\begin{tabular*}{\textwidth}{@{\extracolsep{4in minus 4in}}ld{2.1}d{2.1}d{2.1}
d{2.1}d{2.1}
d{2.1}d{2.1}
d{2.1}d{2.1}
d{2.1}d{2.1}
d{2.1}d{2.1}
d{3.1}d{3.1}
d{3.1}d{3.1}
d{3.1}d{3.1}@{}}
\hline
&& \multicolumn{6}{c}{$\bolds{n=100}$} & \multicolumn{6}{c}{$\bolds{n=250}$}
& \multicolumn{6}{c@{}}{$\bolds{n=500}$} \\[-4pt]
& &\multicolumn{6}{c}{\hrulefill} & \multicolumn{6}{c}{\hrulefill}
& \multicolumn{6}{c@{}}{\hrulefill} \\
& &\multicolumn{2}{c}{$\bolds{{q}_{{n}}}$} &
\multicolumn{2}{c}{$\bolds{{q}_{{  n}}^{{ 0}}}$} & \multicolumn{2}{c}{${\operatorname{\mathbf{GLR}}}$} & \multicolumn{2}{c}{$\bolds{{q}_{{  n}}}$} &
\multicolumn{2}{c}{$%
\bolds{{  q}_{{  n}}^{{ 0}}}$} & \multicolumn{2}{c}{${\operatorname{\mathbf{GLR}}}$} & \multicolumn{2}{c}{$\bolds{{  q}_{{  n}}}$} &
\multicolumn{2}{c}{$\bolds{{  q}_{{  n}}^{{ 0}}}$} & \multicolumn{2}{c@{}}{${\operatorname{\mathbf{GLR}}}$}
\\[-4pt]
& &\multicolumn{2}{c}{\hrulefill}&\multicolumn{2}{c}{\hrulefill}&
\multicolumn{2}{c}{\hrulefill}&
\multicolumn{2}{c}{\hrulefill}&
\multicolumn{2}{c}{\hrulefill}&
\multicolumn{2}{c}{\hrulefill}&
\multicolumn{2}{c}{\hrulefill}&
\multicolumn{2}{c}{\hrulefill}&
\multicolumn{2}{c@{}}{\hrulefill}\\
\multicolumn{1}{@{}l}{$\bolds{(a,\beta)}$} &\multicolumn{1}{c}{$\bolds{\theta}$}& \multicolumn{1}{c}{\textbf{{10\%}}} & \multicolumn{1}{c}{\textbf{{5\%}}} & \multicolumn{1}{c}{\textbf{{10\%}}} &
\multicolumn{1}{c}{\textbf{{5\%}}} & \multicolumn{1}{c}{\textbf{{10\%}}} & \multicolumn{1}{c}{\textbf{{5\%}}} & \multicolumn{1}{c}{\textbf{{10\%}}} &
\multicolumn{1}{c}{\textbf{{5\%}}} &
\multicolumn{1}{c}{\textbf{{10\%}}} & \multicolumn{1}{c}{\textbf{{5\%}}} & \multicolumn{1}{c}{\textbf{{10\%}}} & \multicolumn{1}{c}{\textbf{{5\%}}} &
\multicolumn{1}{c}{\textbf{{10\%}}}
& \multicolumn{1}{c}{\textbf{{5\%}}} & \multicolumn{1}{c}{\textbf{{10\%}}} & \multicolumn{1}{c}{\textbf{{5\%}}} & \multicolumn{1}{c}{\textbf{{10\%}}} &
\multicolumn{1}{c@{}}{\textbf{{5\%}}}
\\
\hline
&&\multicolumn{18}{c}{DGP P.2: Threshold regression} \\
 $(0.0,1.0)$ &  -1.0 &  73.0 &  60.2 &
73.2 &  61.4 &  60.8 &  46.4 &  98.8 &
 96.2 &  98.8 &  96.4 &
95.0 &  90.8 &  100.0 &
100.0
&  100.0 &  100.0 &  100.0 &  100.0 \\
&  -0.5 &  28.0 &  18.0 &  28.4 &
18.4 &  24.2 &  15.6 &  51.8 &  38.2 &
 52.6 &  39.8 &  39.2 &
 28.6 &  81.6 &  69.0 &
81.8 &  69.4 &  62.6 &  51.6 \\
&  -0.2 &  13.4 &  7.2 &  13.4 &
 7.2
&  13.0 &  7.8 &  17.0 &  8.8 &
 16.4
&  8.8 &  13.4 &  6.8 &
24.8 &  15.4 &  25.2 &  15.6 &  18.6 &
 12.2 \\
&  0.2 &  11.0 &  5.8 &  11.4 &
 6.4
&  12.6 &  5.8 &  14.6 &  8.0 &
 15.0
&  8.4 &  12.6 &  7.2 &
23.2 &  14.6 &  23.4 &  14.4 &  18.6 &
 11.0 \\
&  0.5 &  25.0 &  15.6 &  25.2 &
 14.4%
 &  21.0 &  11.8 &  50.4 &  36.8 &
50.4 &  37.0 &  37.2 &
25.4 &
 79.2 &  69.6 &  79.0 &  69.4 &
 64.8
&  51.8 \\
&  1.0 &  71.4 &  59.0 &  72.2 &
 58.6%
 &  60.2 &  43.4 &  97.6 &
 94.6 &  97.6 &  94.6 &  92.6 &
 87.2
&  100.0 &  100.0 &  100.0 &  100.0 &
 100.0 &  100.0 \\ [3pt]
 $(0.2,1.0)$ &  -1.0 &  74.0 &  61.6 &
74.8 &  62.2 &  60.8 &  46.4 &  99.0 &
 96.8 &  98.8 &  96.8 &
95.0 &  90.8 &  100.0 &  100.0 &
100.0 &
 100.0 &  100.0 &  100.0 \\
&  -0.5 &  29.4 &  18.2 &
 29.8 &  18.6 &  24.2 &  15.6 &
 53.0
&  40.2 &  53.6 &  39.8 &
39.2 &  28.6 &  81.0 &  70.4 &  81.8 &
 69.8 &  62.6 &  51.6 \\
&  -0.2 &  14.0 &  7.4 &
13.8 &  7.2 &  13.0 &  7.8 &  17.6 &
 8.4 &  17.2 &  8.2 &
 13.4%
 &  6.8 &  25.2 &  15.2 &  25.2 &
15.6 &  18.6 &  12.2 \\
&  0.2 &  10.8 &  6.2 &
11.2 &  6.0 &  12.6 &  5.8 &  14.0 &
 7.6 &  14.2 &  7.8 &  12.6 &
7.2 &
 23.4 &  14.8 &  23.4 &  14.2 &
 18.6
&  11.0 \\
&  0.5 &  24.0 &  13.8 &
23.8 &  12.8 &  21.0 &  11.8 &  49.0 &
 36.6 &  49.2 &  36.2 &  37.2 &
 25.4
&  78.6 &  68.6 &  78.4 &  69.0 &
64.8 &  51.8 \\
&  1.0 &  70.2 &  56.8 &
70.4 &  56.4 &  60.2 &  43.4 &  97.2 &
 94.2 &  97.4 &  94.6 &  92.6 &
 87.2
&  100.0 &  100.0 &  100.0 &  100.0 &
 100.0 &  100.0 \\[3pt]
 $(0.5,1.0)$ &  -1.0 &  76.2 &  62.2 &
76.8 &  63.0 &  60.8 &  46.4 &  99.0 &
 96.8 &  98.8 &  97.0 &
95.0 &  90.8 &  100.0 &  100.0 &
100.0 &
 100.0 &  100.0 &  100.0 \\
&  -0.5 &  30.4 &  19.2 &
 31.2 &  19.8 &  24.2 &  15.6 &
 55.0
&  41.2 &  56.4 &  40.4 &
39.2 &  28.6 &  81.4 &  71.6 &  81.8 &
 71.4 &  62.6 &  51.6 \\
&  -0.2 &  14.0 &  7.8 &
14.0 &  8.2 &  13.0 &  7.8 &  18.0 &
 9.4 &  18.4 &  10.0 &
13.4 &  6.8 &  26.2 &  16.2 &  25.6 &
 15.8 &  18.6 &  12.2 \\
&  0.2 &  9.4 &  6.0 &
9.6 &  5.8 &  12.6 &  5.8 &  13.0 &
6.8 &  13.8 &  7.0 &  12.6 &  7.2 &
23.0 &  14.0 &  22.4 &  14.0 &  18.6 &
 11.0 \\
&  0.5 &  21.4 &  12.4 &
21.0 &  12.0 &  21.0 &  11.8 &  48.8 &
 35.2 &  48.6 &  35.0 &  37.2 &
 25.4
&  78.4 &  66.8 &  78.0 &  66.8 &
64.8 &  51.8 \\
&  1.0 &  68.8 &  54.6 &
67.8 &  52.8 &  60.2 &  43.4 &  96.8 &
 94.0 &  97.0 &  94.0 &  92.6 &
 87.2
&  100.0 &  100.0 &  100.0 &  100.0 &
 100.0 &  100.0 \\[3pt]
 $(1.0,1.0)$ &  -1.0 &  77.0 &  63.6 &
77.8 &  64.2 &  60.8 &  46.4 &  99.0 &
 96.8 &  99.0 &  96.8 &
95.0 &  90.8 &  100.0 &  100.0 &
100.0 &
 100.0 &  100.0 &  100.0 \\
&  -0.5 &  31.2 &  21.4 &
 32.4 &  22.0 &  24.2 &  15.6 &
 57.6
&  44.0 &  57.8 &  44.0 &
39.2 &  28.6 &  82.6 &  74.2 &  83.2 &
 73.8 &  62.6 &  51.6 \\
&  -0.2 &  14.0 &  8.0 &
14.8 &  8.2 &  13.0 &  7.8
&
 18.6 &  11.0 &  19.6 &
10.8 &  13.4 &  6.8 &  26.0 &  16.0 &
 26.6 &  17.0 &  18.6 &  12.2 \\
&  0.2 &  9.2 &  6.0 &
9.4 &  6.0 &  12.6 &  5.8 &  12.2 &
6.4 &  12.4 &  6.8 &  12.6 &  7.2 &
21.4 &  12.6 &  21.0 &  13.0 &  18.6 &
 11.0 \\
&  0.5 &  19.2 &  10.0 &
18.4 &  9.6 &  21.0 &  11.8 &  47.2 &
 31.6 &  46.8 &  32.2 &  37.2 &
 25.4
&  76.4 &  64.0 &  75.8 &  65.2 &
64.8 &  51.8 \\
&  1.0 &  63.8 &  48.2 &
62.2 &  46.2 &  60.2 &  43.4 &  95.8 &
 92.8 &  96.0 &  93.0 &  92.6 &
 87.2
&  100.0 &  100.0 &  100.0 &  100.0 &
 100.0 &  100.0 \\
 \hline
\end{tabular*}}}}
\tabnotetext[]{}{
{\fontsize{7.2}{9.2}{\selectfont{
\textit{Notes}:
(i) 1000 iterations;
(ii) $\operatorname{GLR}$ the generalized likelihood ratio
test, ${  q}_{{  n}}$ and ${  q}_{{
n}}^{%
{ 0}}$ loss function-based tests;
(iii) ${  q}_{{  n}}$ is standardized by
$\operatorname{SSR}_{1}$, the sum of squared residuals of the
nonparametric regression estimates, and ${  q}_{{  n}}^{{
 0}}$ is standardized by ${\operatorname{SSR}}_{0}$, the sum of
squared residuals of the null linear model;
(iv) The uniform kernel is used for ${\operatorname{GLR}}$,
${q}_{{  n}}$ and ${  q}_{{  n}}^{{ 0}}$; the bandwidth $  h=S_{{  X}}{  n}^{-2/9}$,
where ${  S}_{{  X}}$ is the sample standard
deviation
of $\{{  X}_{t}\}_{t=1}^{n};$
(v) The ${  q}_{{  n}}$ tests are based on
the
linex loss function: ${d(z)=\frac{\beta}{\alpha^{2}}[\exp(
\alpha
z)-1-\alpha z];}$
(vi) DGP P.2, ${Y}_{t} =1+X_{t}
1(X_{t}%
 >0)+(1+\theta)X_{t} 1(X_{t} \leq0)+
\varepsilon_{t}$, where $ \{ \varepsilon_{i} \}
\sim
\mathrm{i.i.d.}\ N(0,1).$}}}}
\end{sidewaystable}

\begin{sidewaystable}
\tabcolsep=0pt
\tablewidth=\textwidth
\caption{Empirical powers of tests}\label{tab6}
{\fontsize{8.2}{10.2}{\selectfont{\begin{tabular*}{\textwidth}{@{\extracolsep{\fill}}ld{2.1}
d{2.1}d{2.1}
d{2.1}d{2.1}
d{2.1}d{2.1}
d{2.1}d{2.1}
d{2.1}d{2.1}
d{2.1}d{2.1}
d{3.1}d{3.1}
d{3.1}d{3.1}
d{3.1}d{3.1}@{}}
\hline
&& \multicolumn{6}{c}{$\bolds{n=100}$} & \multicolumn{6}{c}{$\bolds{n=250}$}
& \multicolumn{6}{c@{}}{$\bolds{n=500}$} \\[-4pt]
& &\multicolumn{6}{c}{\hrulefill} & \multicolumn{6}{c}{\hrulefill}
& \multicolumn{6}{c@{}}{\hrulefill} \\
& &\multicolumn{2}{c}{$\bolds{{q}_{{n}}}$} &
\multicolumn{2}{c}{$\bolds{{q}_{{  n}}^{{ 0}}}$} & \multicolumn{2}{c}{${\operatorname{\mathbf{GLR}}}$} & \multicolumn{2}{c}{$\bolds{{q}_{{  n}}}$} &
\multicolumn{2}{c}{$%
\bolds{{  q}_{{  n}}^{{ 0}}}$} & \multicolumn{2}{c}{${\operatorname{\mathbf{GLR}}}$} & \multicolumn{2}{c}{$\bolds{{  q}_{{  n}}}$} &
\multicolumn{2}{c}{$\bolds{{  q}_{{  n}}^{{ 0}}}$} & \multicolumn{2}{c@{}}{${\operatorname{\mathbf{GLR}}}$}
\\[-4pt]
& &\multicolumn{2}{c}{\hrulefill}&\multicolumn{2}{c}{\hrulefill}&
\multicolumn{2}{c}{\hrulefill}&
\multicolumn{2}{c}{\hrulefill}&
\multicolumn{2}{c}{\hrulefill}&
\multicolumn{2}{c}{\hrulefill}&
\multicolumn{2}{c}{\hrulefill}&
\multicolumn{2}{c}{\hrulefill}&
\multicolumn{2}{c@{}}{\hrulefill}\\
\multicolumn{1}{@{}l}{$\bolds{(a,\beta)}$} &\multicolumn{1}{c}{$\bolds{\theta}$}& \multicolumn{1}{c}{\textbf{{10\%}}} & \multicolumn{1}{c}{\textbf{{5\%}}} & \multicolumn{1}{c}{\textbf{{10\%}}} &
\multicolumn{1}{c}{\textbf{{5\%}}} & \multicolumn{1}{c}{\textbf{{10\%}}} & \multicolumn{1}{c}{\textbf{{5\%}}} & \multicolumn{1}{c}{\textbf{{10\%}}} &
\multicolumn{1}{c}{\textbf{{5\%}}} &
\multicolumn{1}{c}{\textbf{{10\%}}} & \multicolumn{1}{c}{\textbf{{5\%}}} & \multicolumn{1}{c}{\textbf{{10\%}}} & \multicolumn{1}{c}{\textbf{{5\%}}} &
\multicolumn{1}{c}{\textbf{{10\%}}}
& \multicolumn{1}{c}{\textbf{{5\%}}} & \multicolumn{1}{c}{\textbf{{10\%}}} & \multicolumn{1}{c}{\textbf{{5\%}}} & \multicolumn{1}{c}{\textbf{{10\%}}} &
\multicolumn{1}{c@{}}{\textbf{{5\%}}}
\\
\hline
&&\multicolumn{18}{c}{DGP P.3: Smooth transition regression}\\
 $(0.0,1.0)$ &  -1.0 &  53.2 &  40.0 &
54.0 &  41.6 &  43.8 &  33.6 &  90.2 &
 84.6 &  90.2 &  84.2 &  78.8 &
 69.8
&  99.6 &  99.2 &  99.6 &  99.2 &
97.8 &  95.8 \\
&  -0.5 &  22.6 &  13.0 &  22.8 &
13.8 &  20.8 &  11.8 &  37.8 &  27.2 &
 38.8 &  27.4 &  27.8 &  20.4 &
 61.2
&  49.2 &  60.2 &  49.8 &  47.6 &
35.2 \\
&  0.5 &  20.6 &  10.4 &  19.8 &
 11.4%
 &  16.6 &  9.6 &  36.2 &  25.0 &
35.6 &  26.2 &  29.2 &  17.0 &  60.8 &
 49.4 &  60.2 &  50.6 &  46.8 &
 35.0
\\
&  1.0 &  52.2 &  38.2 &  51.4 &
 37.6%
 &  44.0 &  29.8 &  87.4 &  79.6 &
86.8 &  79.8 &  78.2 &  69.0 &  99.4 &
 98.6 &  99.4 &  98.6 &  98.4 &
 96.0
\\
&  1.5 &  85.6 &  75.2 &  86.6 &
 77.0%
 &  75.6 &  65.6 &  99.6 &  99.4 &
99.6 &  99.2 &  98.2 &  96.8 &  100.0
&
 100.0 &  100.0 &  100.0 &  100.0 &
100.0 \\[3pt]
 $(0.2,1.0)$ &  -1.0 &  55.0 &  40.4 &
56.0 &  42.6 &  43.8 &  33.6 &  90.6 &
 85.0 &  90.4 &  85.6 &  78.8 &
 69.8
&  99.6 &  99.2 &  99.6 &  99.2 &
97.8 &  95.8 \\
&  -0.5 &  23.8 &  14.0 &  23.0 &
14.2 &  20.8 &  11.8 &  39.4 &  28.4 &
 40.2 &  28.2 &  27.8 &  20.4 &
 61.6
&  49.8 &  60.8 &  51.0 &  47.6 &
35.2 \\
&  0.5 &  19.4 &  10.2 &  19.4 &
 10.4%
 &  16.6 &  9.6 &  35.4 &  24.4 &
35.2 &  25.2 &  29.2 &  17.0 &  59.2 &
 48.8 &  59.0 &  49.2 &  46.8 &
 35.0
\\
&  1.0 &  50.4 &  37.2 &  49.8 &
 36.2%
 &  44.0 &  29.8 &  86.8 &  79.0 &
85.6 &  78.4 &  78.2 &  69.0 &  99.4 &
 98.6 &  99.4 &  98.6 &  98.4 &
 96.0
\\
&  1.5 &  84.6 &  74.2 &  85.0 &
 75.2%
 &  75.6 &  65.6 &  99.4 &  99.4 &
99.6 &  99.0 &  98.2 &  96.8 &  100.0
&
 100.0 &  100.0 &  100.0 &  100.0 &
100.0 \\[3pt]
 $(0.5,1.0)$ &  -1.0 &  56.6 &  43.0 &
57.6 &  44.0 &  43.8 &  33.6 &  86.2 &
 63.6 &  86.0 &  64.6 &  78.8 &
 69.8
&  99.6 &  99.2 &  99.6 &  99.2 &
97.8 &  95.8 \\
&  -0.5 &  25.0 &  14.4 &  24.0 &
14.4 &  20.8 &  11.8 &  40.4 &  29.4 &
 40.6 &  29.2 &  27.8 &  20.4 &
 62.4
&  51.4 &  61.4 &  51.8 &  47.6 &
35.2 \\
&  0.5 &  18.4 &  9.4 &  18.4 &
 8.8
&  16.6 &  9.6 &  34.6 &  23.0 &
 34.0%
 &  23.2 &  29.2 &  17.0 &  58.0 &
47.4 &  57.8 &  48.0 &  46.8 &  35.0 \\
&  1.0 &  48.8 &  34.0 &  47.6 &
 33.2%
 &  44.0 &  29.8 &  85.0 &  77.6 &
85.0 &  77.8 &  78.2 &  69.0 &  99.2 &
 98.2 &  99.4 &  98.4 &  98.4 &
 96.0
\\
&  1.5 &  83.0 &  70.2 &  82.6 &
 72.0%
 &  75.6 &  65.6 &  99.4 &  98.8 &
99.6 &  99.0 &  98.2 &  96.8 &  100.0
&
 100.0 &  100.0 &  100.0 &  100.0 &
100.0 \\[3pt]
 $(1.0,1.0)$ &  -1.0 &  58.4 &  46.4 &
58.8 &  46.0 &  43.8 &  33.6 &  92.0 &
 87.4 &  91.8 &  86.8 &  78.8 &
 69.8
&  99.6 &  99.2 &  99.6 &  99.2 &
97.8 &  95.8 \\
&  -0.5 &  26.2 &  16.2 &  26.2 &
16.4 &  20.8 &  11.8 &  43.6 &  30.8 &
 43.8 &  30.6 &  27.8 &  20.4 &
 64.6
&  52.4 &  64.0 &  52.2 &  47.6 &
35.2 \\
&  0.5 &  14.6 &  7.4 &  15.0 &
 7.6
&  16.6 &  9.6 &  31.0 &  20.2 &
 31.8%
 &  20.6 &  29.2 &  17.0 &  56.8 &
44.4 &  56.6 &  45.2 &  46.8 &  35.0 \\
&  1.0 &  44.0 &  27.8 &  43.2 &
 29.0%
 &  44.0 &  29.8 &  83.0 &  75.4 &
82.8 &  75.8 &  78.2 &  69.0 &  99.2 &
 97.6 &  99.2 &  97.8 &  98.4 &
 96.0
\\
&  1.5 &  78.2 &  62.8 &  79.4 &
 64.8%
 &  75.6 &  65.6 &  99.4 &  98.4 &
99.4 &  98.6 &  98.2 &  96.8 &  100.0
&
 100.0 &  100.0 &  100.0 &  100.0 &
100.0 \\
\hline
\end{tabular*}}}}
\tabnotetext[]{}{
{\fontsize{8}{10}{\selectfont{
\textit{Notes}:
(i) 1000 iterations;
(ii) $\operatorname{GLR}$ the generalized likelihood ratio
test, ${q}_{{  n}}$ and ${  q}_{{
n}}^{%
{ 0}}$, loss function-based tests;
(iii) ${  q}_{{  n}}$ is standardized by
${\operatorname{SSR}}_{1}$, the sum of squared residuals of the
nonparametric regression estimates, and ${ q}_{{  n}}^{{
 0}}$ is standardized by ${\operatorname{SSR}}_{0}$, the sum of
squared residuals of the null linear model;
(iv) The uniform kernel is used for ${\operatorname{GLR}}$, ${q}_{{ n}}$ and
${q}_{{  n}}^{{ 0}}$;
the bandwidth $h=S_{{X}}{  n}^{-2/9}$, where ${  S}_{{  X}}$ is the sample standard
deviation
of $\{  X_{t}\}_{t=1}^{n}$;
(v) The ${  q}_{{  n}}$ tests are based on
the
linex loss function: $d(z)=\frac{\beta}{\alpha^{2}}[\exp(
\alpha
z)-1-\alpha z];$
(vi) DGP P.3, $Y_{t}=1+X_{t}+
[1-F(X_{t})\theta]X_{t}+\varepsilon_{t}$, $F(X_{t})=
\frac{1}{ 1+\exp (-X_{t})}$, where $\{\varepsilon_{i} \} \sim
\mathrm{i.i.d.}\ N(0,1)$.}}}}
\end{sidewaystable}

Next, we examine the size of the tests based on the bootstrap. Table~\ref{tab3}
shows that overall, the rejection rates of all tests based on the bootstrap
are close to the significance levels (10\% and 5\%), indicating the
gain of
using the bootstrap in finite samples. The sizes of all tests are
robust to
a variety of error distributions, confirming the Wilks phenomena that the
asymptotic distribution of both the $q_{n}$ and $\lambda_{n}$ tests are
distribution free. For the loss function tests $q_{n}$ and $q_{n}^{0}$, the
sizes are very similar for different choices of parameters $(\alpha,\beta)$
governing the shape of the linex loss function. We note that when asymptotic
critical values are used, the sizes of the tests with $h=S_{X}n^{-2/9}$ are
slightly better than with $h=S_{X}n^{-1/5}$. When bootstrap critical values
are used, however, the sizes of all tests with $h=S_{X}n^{-2/9}$ and $%
h=S_{X}n^{-1/5}$, respectively, are very similar.

Next, we turn to the powers of the tests under $\mathbb{H}_{A}$. Since the
sizes of the tests using asymptotic critical values are different in finite
samples, we use the bootstrap procedure only, which delivers similar sizes
close to significance levels and thus provides a fair ground for comparison.
Tables~\ref{tab4}--\ref{tab6} report the empirical rejection rates of the tests under
DGP~\ref{dgp1}
(quadratic regression), DGP~\ref{dgp2} (threshold regression) and DGP~\ref{dgp3} (smooth
transition regression), respectively. For all DGPs, the loss function
tests $%
q_{n}$ and $q_{n}^{0}$ are more powerful than the GLR test, confirming our
asymptotic efficiency analysis. Interestingly, for the two loss function
tests, $q_{n}$, which is standardized by the nonparametric $\operatorname{SSR}_{1}$, is
roughly equally powerful to $q_{n}^{0}$, which is standardized by the
parametric $\operatorname{SSR}_{0}$, although asymptotic analysis suggests that $q_{n}$
should be more powerful than $q_{n}^{0}$ under $\mathbb{H}_{A}$,
because $%
\operatorname{SSR}_{0}$ is significantly larger than $\operatorname{SSR}_{1}$ under $\mathbb{H}_{A}$.
Obviously, this is due to the use of the bootstrap. Since the bootstrap
statistics $q_{n}^{\ast}$ and $q_{n}^{0\ast}$ are standardized by $%
\operatorname{SSR}_{0}^{\ast}$ and $\operatorname{SSR}_{1}^{\ast}$, respectively, where
$\operatorname{SSR}_{1}^{\ast
}<\operatorname{SSR}_{0}^{\ast}$, the ranking between $q_{n}$ and $q_{n}^{\ast}$ remains
more or less similar to the ranking between $q_{n}^{0}$ and
$q_{n}^{0\ast}$, and therefore $q_{n}^{\ast}$ and $q_{n}^{0\ast}$ have similar power.
Under DGP~\ref{dgp1}, the powers of $q_{n}$ and $q_{n}^{0}$ increase as the
degree of
asymmetry of the linex loss function, which is indexed by $\alpha$,
increases. When $\alpha=1$, the powers of $q_{n}$ and $q_{n}^{0}$ are
substantially higher than the GLR test. Under DGP~\ref{dgp2}, there is some tendency
that the powers of $q_{n}$ and $q_{n}^{0}$ increase in $\alpha$ for
$\theta
<0$, whereas they decrease in $\alpha$ for $\theta>0$. When $\theta$ is
close to 0, $q_{n}$ and $q_{n}^{0}$ have similar power to the GLR test, but
as $|\theta|>0$ increases, they are more powerful than the GLR test.
Similarly, under DGP~\ref{dgp3}, the powers of $q_{n}$ and $q_{n}^{0}$ increase
in $%
\alpha$ for $\theta<0$, whereas they decrease in $\alpha$ for
$\theta>0$.
Nevertheless, by and large, the powers of both $q_{n}$ and $q_{n}^{0}$
do not change much across the
different choices of parameters $(\alpha,\beta)$
governing the shape of the linex loss function. Although the shape of the
loss function changes dramatically when $\alpha$ changes from 0 to 1, the
powers of $q_{n}$ and $q_{n}^{0}$ remain relatively robust.

All tests become more powerful as the departures from linearity increases
(as characterized by the value of $\theta$ in each DGP), and as the sample
size $n$ increases.

\section{Conclusion}\label{sec7}
The GLR test has been proposed as a generally applicable method for
nonparametric testing problems. It inherits many advantages of the
maximum LR test for parametric models. In this paper, we have shown
that despite its general nature and many appealing features, the GLR
test does not have the optimal power property of the classical LR test.
We propose a loss function test in a time series context. The new test
enjoys the same appealing features as the GLR test, but is more
powerful in terms of Pitman's asymptotic efficiency. This holds no
matter what kernel and bandwidth are used, and even when the true
likelihood function is available for the GLR test. The efficiency gain,
together with more relevance to decision making under uncertainty of
using a loss function, suggests that the loss function approach can be
a generally applicable and
powerful nonparametric inference procedure alternative to the GLR
principle.%

\begin{appendix}\label{app}
\section*{Mathematical appendix}

Throughout the appendix, we let
$\tilde{m%
}_{h}(x)$ be defined in the same way as $\hat{m}_{h}(x)$ in (\ref{eq4.2}) with
$\{
\varepsilon_{t}=Y_{t}-g_{0}(X_{t})\}_{t=1}^{n}$ replacing $\{ \hat{%
\varepsilon}_{t}=Y_{t}-g(X_{t},\hat{\theta}_{0})\}_{t=1}^{n}$. Also,
$C\in
(1,\infty)$ denotes a generic bounded constant. This appendix provides the
structure of our proof strategy. We leave the detailed proofs of most
technical lemmas and propositions to the supplementary material.

\begin{pf*}{Proof of Theorem~\ref{th1}}
Theorem~\ref{th1} follows as a special case
of Theorem~\ref{th3}(i) with $\delta(X_{t})=0$.
\end{pf*}
\begin{pf*}{Proof of Theorem~\ref{th2}}
Theorem~\ref{th2} follows as a special case
of Theorem~\ref{th3}(ii) with $\delta(X_{t})=0$.
\end{pf*}

\begin{pf*}{Proof of Theorem~\ref{th3}\normalfont{(i)}} We shall first derive the
asymptotic distribution of $q_{n}$ under $\mathbb{H}_{n}(a_{n})$. From
Lemmas~\ref{le3.1} and~\ref{le3.2} and Propositions~\ref{pr3.1} and~\ref{pr3.2} below, we can obtain
\[
\frac{h^{{p}/{2}}D^{-1}q_{n}-h^{-{p}/{2}}\sigma^{-2}\int
K^{2}(u)\,du\int\sigma^{2}(x)\,dx}{\sqrt{2\sigma^{-4}\int [ \int
K(u)K(u+v)\,du ] ^{2}\,dv\int\sigma^{4}(x)\,dx}}\stackrel{d}
{\longrightarrow}N(\psi,1).
\]
The desired result of Theorem~\ref{th3}(i) then follows immediately.
\end{pf*}

\renewcommand{\thelemma}{\Alph{section}.\arabic{lemma}}
\begin{lemma}\label{le3.1}
Under the conditions of Theorem~\ref{th3},
$\hat{Q}_{n}=D\sum_{t=1}^{n}\hat{m}_{h}^{2}(X_{t})+o_{p}(h^{-p/2})$.
\end{lemma}

\begin{lemma}\label{le3.2}
Under the conditions of Theorem~\ref{th3}, $\sum_{t=1}^{n}\hat{m}_{h}^{2}(X_{t})=n\times \int\hat{m}%
_{h}^{2}(x)f(x)\,dx+o_{p}(h^{-p/2})$.
\end{lemma}

\renewcommand{\theproposition}{\Alph{section}.\arabic{proposition}}
\begin{proposition}\label{pr3.1}
$\!\!\!$Under the conditions of Theorem
\ref{th3}, $n\int\hat{m}_{h}^{2}(x)f(x)\,dx\!=n\int\tilde{m}%
_{h}^{2}(x)f(x)\,dx+h^{-p/2}E[\delta^{2}(X_{t})]+o_{p}(h^{-p/2})$.
\end{proposition}

\begin{proposition}\label{pr3.2}
Under the conditions of Theorem
\ref{th3}, and $\mathbb{H}_{n}(a_{n})$ with $a_{n}=n^{-1/2}h^{-p/4}$,%
\begin{eqnarray*}
&&\biggl[ nh^{p/2}\int\tilde{m}_{h}^{2}(x)f(x)\,dx-h^{-{p}/{2}}a(K)
\int\sigma^{2}(x)\,dx \biggr] \Big/\sqrt{2b(K)\int\sigma
^{4}(x)\,dx}\\
&&\qquad\stackrel{d} {%
\longrightarrow}N(\psi,1).
\end{eqnarray*}
\end{proposition}

\begin{pf*}{Proof of Lemma~\ref{le3.1}}
Given in the supplementary material.
\end{pf*}

\begin{pf*}{Proof of Lemma~\ref{le3.2}} Given in the supplementary material.
\end{pf*}

\begin{pf*}{Proof of Proposition~\ref{pr3.1}} Given in the supplementary
material.~
\end{pf*}

\begin{pf*}{Proof of Proposition~\ref{pr3.2}} Proposition~\ref{pr3.2} follows from
Lemmas~\ref{le3.3} and~\ref{le3.4} below.
\end{pf*}

\begin{lemma}\label{le3.3}
Put $\hat{H}_{q}=n^{-1}\sum_{t=2}^{n}%
\sum_{s=1}^{t-1}H_{n}(Z_{t},Z_{s})$, where $Z_{t}=(\varepsilon
_{t},X_{t}^{\prime})^{\prime}$, $H_{n}(Z_{t},Z_{s})=2\varepsilon
_{t}\varepsilon_{s}W_{h}(X_{t},X_{s})$, and
\[
W_{h}(X_{t},X_{s})=\int
\frac{\mathbf{K}_{h}(X_{t}-x)\mathbf{K}_{h}(X_{s}-x)%
}{f(x)}\,dx.
\]
Suppose Assumptions~\ref{assa1} and~\ref{assa4} hold, $h\propto n^{-\omega}$
for $\omega\in(0,1/2p)$ and $p<4$. Then
\[
n\int\hat{m}_{h}^{2}(x)g(x)\,dx=h^{-p}\int
\mathbf{K}^{2}(\mathbf{u})\,d\mathbf{%
u}\int\sigma
^{2}(x)\,dx+\hat{H}_{q}+o_{p}\bigl(h^{-p/2}
\bigr).
\]
\end{lemma}

\begin{lemma}\label{le3.4}
Suppose Assumptions~\ref{assa1} and~\ref{assa4} hold,
and $h\propto n^{-\omega}$ for $\omega\in(0,1/2p)$. Define
$V_{q}=2\int [ \int\mathbf{K}(\mathbf{v})\mathbf{K}(\mathbf{u}+\mathbf
{v%
})\,d\mathbf{v} ] ^{2}\,du\int\sigma^{4}(x)\,dx$. Then $V_{q}^{-1/2}\times h^{p/2}\hat{H}_{q}\stackrel{d}{\longrightarrow}N(\psi,1)$.
\end{lemma}

Since Lemmas~\ref{le3.3} and~\ref{le3.4} are the key results for deriving the asymptotic
distributions of the proposed $q_{n}$ statistic when $\{ \varepsilon
_{t}\}$
may not be an i.i.d. sequence nor martingale difference sequence,
we provide detailed proofs for them below.\vadjust{\goodbreak}

\begin{pf*}{Proof of Lemma~\ref{le3.3}} Let $\hat{F}_{n}(x)$ be the empirical
distribution function of $\{X_{t}\}_{t=1}^{n}$. We have
\begin{eqnarray}\label{eqA.1}
&& n\int\tilde{m}_{h}^{2}(x)f(x)\,dx\nonumber\\
&&\qquad
=n\int\frac{ [ n^{-1}\sum_{s=1}^{n}\varepsilon_{s}\mathbf{K}%
_{h}(X_{t}-X_{s}) ] ^{2}}{f(x)}\,dx\nonumber\\
&&\qquad\quad{}+\int\Biggl[ n^{-1}\sum
_{s=1}^{n}%
\varepsilon_{s}
\mathbf{K}_{h}(X_{t}-X_{s}) \Biggr]
^{2} \biggl[ \frac{1}{\hat{f%
}^{2}(x)}-\frac{1}{f^{2}(x)} \biggr] f(x)\,dx
\nonumber\\
&&\qquad=n^{-1}\sum_{t=1}^{n}\sum
_{s=1}^{n}\varepsilon_{t}
\varepsilon_{s}\int\frac{\mathbf{K}_{h}(X_{t}-x)\mathbf
{K}_{h}(X_{s}-x)}{f(x)}%
\,dx
\nonumber
\\[-8pt]
\\[-8pt]
\nonumber
&&\qquad\quad{}+O_{p}
\bigl(n^{-1}h^{-p}\bigr)O_{p}\bigl(n^{-1/2}h^{-p/2}
\ln n+h^{2}\bigr)
\\
&&\qquad=n^{-1}\sum_{t=1}^{n}
\varepsilon_{t}^{2}\int\frac{\mathbf{K}%
_{h}^{2}(X_{t}-x)}{f(x)}\,dx\nonumber\\
&&\qquad\quad{}+n^{-1}
\sum_{1\leq s\leq t\leq n}^{n}2\varepsilon_{t}
\varepsilon_{s}\int\frac{\mathbf{K}_{h}(X_{t}-x)\mathbf
{K}_{h}(X_{s}-x)%
}{f(x)}+o_{p}
\bigl(h^{-p/2}\bigr)
\nonumber
\\
&&\qquad=\hat{C}_{q}+\hat{H}_{q}+o_{P}
\bigl(h^{-p/2}\bigr), \nonumber
\end{eqnarray}
where we have made use of the fact that $\sup_{x\in\mathbb{G}}\llvert
\hat{f}(x)-f(x)\rrvert =O_{p}(n^{-1/2}\times  h^{-p/2}\ln n+h^{2})$ given
Assumption~\ref{assa2}, and $h\propto n^{-\omega}$ for $\omega\in(0,1/2p)$.

By change of variable, the law of iterated expectations, and Assumption~\ref{assa1},
we can obtain%
%
\begin{eqnarray}\label{eqA.2}
E(\hat{C}_{q})&=&\iint\sigma^{2}(x)\frac{\mathbf{K}_{h}^{2}(y-x)}{f(x)}%
f(y)\,dx\,dy
\nonumber
\\[-8pt]
\\[-8pt]
\nonumber
&=&h^{-p}\int\mathbf{K}^{2}(u)\,du\int\sigma
^{2}(x)\,dx \bigl[ 1+O\bigl(h^{2}\bigr)%
\bigr].
\end{eqnarray}
On the other hand, by Chebyshev's inequality and the fact that $E(\hat
{C}%
_{q}-E\hat{C}_{q})^{2}=O_{p}(n^{-1}h^{-2p})$ given Assumption~\ref{assa1}, we
have%
%
\begin{equation}
\hat{C}_{q}=E(\hat{C}_{q})+O_{p}
\bigl(n^{-1/2}h^{-p}\bigr). \label{eqA.3}
\end{equation}
Combining (\ref{eqA.1})--(\ref{eqA.3}) and $p<4$ then yields the desired result of
Lemma~\ref{le3.3}.
\end{pf*}%

\begin{pf*}{Proof of Lemma~\ref{le3.4}} Because $%
E[H_{n}(Z_{t},z)]=E[H_{n}(z^{\prime},Z_{s})]=0$ for all $z,z^{\prime
},\hat{%
H}_{q}\equiv n^{-1}\sum_{1\leq s<t\leq n}H_{n}(Z_{t},Z_{s})$ is a degenerate
$U$-statistic. Following Tenreiro's (\citeyear{Ten97}) central limit\vadjust{\goodbreak}
theorem for
degenerate $U$-statistics of a time series context process, we have $%
[n^{-2}\sum_{1\leq s<t\leq
n}E[h^{p}H_{n}^{2}(Z_{t},Z_{s})]]^{-1/2}h^{p/2}%
\hat{H}_{q}\stackrel{d}{\rightarrow}N(0,1)$ as $n\rightarrow\infty$
if the
following conditions are satisfied: For some constants $\delta
_{0}>0,\gamma
_{0}<\frac{1}{2}$ and $\gamma_{1}>0$, (i)~$u_{n}(4+\delta
_{0})=O(n^{\gamma
_{0}})$, (ii)~$v_{n}(2)=o(1)$, (iii)~$w_{n}(2+\frac{\delta
_{0}}{2})=o(n^{%
{1}/{2}})$ and (iv)~$z_{n}(2)n^{\gamma_{1}}=O(1)$, where
\begin{eqnarray*}
u_{n}(r) &=&h^{p/2}\max\Bigl\{ \max_{1\leq t\leq
n}
\bigl\Vert H_{n}(Z_{t},Z_{0})\bigr\Vert_{r},
\bigl\Vert H_{n}(Z_{0},\bar{Z}_{0})\bigr\Vert
_{r} \Bigr\},
\\
v_{n}(r) &=&h^{p}\max\Bigl\{ \max_{1\leq t\leq
n}
\bigl\Vert G_{n0}(Z_{t},Z_{0})\bigr\Vert_{r},
\bigl\Vert G_{n0}(Z_{0},\bar{Z}_{0})\bigr\Vert
_{r} \Bigr\},
\\
w_{n}(r) &=&h^{p}\bigl\Vert G_{n0}(Z_{0},Z_{0})
\bigr\Vert_{r},
\\
z_{n}(r) &=&h^{p}\max_{0\leq t\leq n,1\leq s\leq n}\max\bigl\{
\bigl\Vert G_{ns}(Z_{t},Z_{0})\bigr\Vert_{r},
\bigl\Vert G_{ns}(Z_{0},Z_{t})\bigr\Vert_{r},
\bigl\Vert G_{ns}(Z_{0},\bar{Z%
}_{0})\bigr\Vert
_{r} \bigr\},
\end{eqnarray*}
$G_{ns}(u,v)=E [ H_{n}(Z_{s},u)H_{n}(Z_{0},v) ] $ for $s\in\mathbb{%
N}$ and $u,v\in\mathbb{R}^{p}$, $\bar{Z}_{0}$ is an independent copy
of $%
Z_{0}$, and $\Vert\xi\Vert_{r}=E^{{1}/{r}}|\xi|^{r}$.

We first show $n^{-2}\sum_{1\leq s<t\leq
n}h^{p}E[H_{n}^{2}(Z_{t},Z_{s})]\rightarrow V_{q}$\vspace*{2pt} as $n\rightarrow
\infty$.
By change of variables and Assumption~\ref{assa1}, it is straightforward to
calculate%
%
\begin{eqnarray}\label{eqA.4}
&&n^{-2}\sum_{1\leq s<t\leq n}h^{p}E
\bigl[H_{n}^{2}(Z_{t},Z_{s})\bigr]\nonumber\\
&&\qquad=4h^{p}n^{-2}\sum_{1\leq s<t\leq n}E
\bigl[\varepsilon_{t}^{2}\varepsilon_{s}^{2}W_{h}^{2}(X_{t},X_{s})
\bigr]
\\
&&\qquad\rightarrow 2\int\biggl[ \int\mathbf{K}(v)\mathbf{K}(u+v)\,dv \biggr]
^{2}\,du\int\sigma^{4}(x)\,dx\equiv V_{q}.\nonumber
\end{eqnarray}

We now verify conditions (i)--(iv). We first consider condition (i). By the
Cauchy--Schwarz inequality and change of variables, we have for all
$t\geq0$,

\begin{eqnarray*}
E\bigl\llvert h^{p/2}H_{n}(Z_{t},Z_{0})
\bigr\rrvert^{r} &=&2^{\gamma}h^{%
({p}/{2})r}E\bigl\llvert
\varepsilon_{t}^{r}\varepsilon_{0}^{r}W_{h}^{r}(X_{t},X_{0})
\bigr\rrvert
\\
&\leq&2^{\gamma}h^{({p}/{2})r}\bigl(E\varepsilon_{0}^{2cr}
\bigr)^{{1}/{c}%
} \bigl( E\bigl\llvert W_{h}(X_{t},X_{0})
\bigr\rrvert^{cr} \bigr) ^{{1}/{c%
}}
\\
&\leq&Ch^{({p}/{2})r} \biggl[ \int\bigl\llvert W_{h}(x,x_{0})
\bigr\rrvert^{cr}f_{X_{t},X_{0}}(x,x_{0})\,dx\,dx_{0}
\biggr] ^{{1}/{c}}
\\
&\leq&Ch^{({p}/{2})r} \bigl( h^{-pcr}h^{p} \bigr)
^{{1}/{c}}\leq Ch^{-%
({r}/{2})p+({p}/{c})}
\end{eqnarray*}
for all $c>1$, and given $E(\varepsilon_{t}^{8+\delta})\leq C$. We
obtain $%
\Vert h^{p/2}H_{n}(Z_{t},Z_{0})\Vert_{r}=\break  ( Ch^{-({r}/{2})p+({p}/{c})%
} ) ^{{1}/{r}}\leq Ch^{-{p}/{2}+{p}/{(cr)}}$. Given $h\propto
n^{-\omega}$ for $\omega\in(0,1/2p)$, we have $%
\Vert h^{p/2}H_{n}(Z_{t},Z_{0})\Vert_{r}\leq Cn^{\omega p(
{1}/{2}-{2}/{(8+\delta)})}$, with $c=\frac{8+\delta}{2r}$ and if $r<4+\frac{\delta
}{2}$.%
 By a similar argument and replacing $f_{X_{t},X_{0}}(x,x_{0})$ with $%
f(x)f(x_{0})$, we can obtain the same order of magnitude for $%
\Vert h^{p/2}H_{n}(Z_{0},\bar{Z}_{0})\Vert_{r}$. Hence, we obtain
$u_{n}(r)\leq
Cn^{\omega p({1}/{2}-{2}/{(8+\delta)})}$, and condition (i)
holds by
setting \mbox{$\gamma_{0}=\omega p(\frac{1}{2}-\frac{2}{8+\delta})$}.

Now we verify condition (ii). Note that for all $s\geq0$, we have%
\begin{eqnarray*}
G_{ns}\bigl(z,z^{\prime}\bigr) &=&E \bigl[ H_{n}(Z_{s},z)H_{n}
\bigl(Z_{0},z^{\prime}\bigr)%
\bigr]
\\
&=&4E \bigl[ \varepsilon_{t}\varepsilon W_{h}(X_{s},x)
\varepsilon_{0}\varepsilon^{\prime}W_{h}
\bigl(X_{0},x^{\prime}\bigr) \bigr]
\\
&=&4\varepsilon\cdot\varepsilon^{\prime}E \bigl[ \varepsilon
_{s}\varepsilon_{0}W_{h}(X_{s},x)W_{h}
\bigl(X_{0},x^{\prime}\bigr) \bigr],
\end{eqnarray*}
where $z=(\varepsilon,x)$ and $z^{\prime}=(\varepsilon^{\prime
},x^{\prime})$. To compute the order of magnitude for $v_{n}(r)$, we first
consider the case of $s=0$. We have%
\begin{eqnarray*}
G_{n0}\bigl(z,z^{\prime}\bigr) &=&4\varepsilon\varepsilon
^{\prime}E_{0} \bigl[ \bar{%
\varepsilon}_{0}^{2}W_{h}(
\bar{X}_{0},x)W_{h}\bigl(\bar{X}_{0},x^{\prime}
\bigr) \bigr]
\\
&=&4\varepsilon\varepsilon^{\prime}E_{0} \bigl[ \sigma
^{2}(\bar{X}%
_{0})W_{h}(
\bar{X}_{0},x)W_{h}\bigl(\bar{X}_{0},x^{\prime}
\bigr) \bigr],
\end{eqnarray*}
where $E_{0}(\cdot)$ is an expectation taken over $(\bar{X}_{0},\bar{%
\varepsilon}_{0})$. By the Cauchy--Schwarz inequality and change of
variables, we have%
\begin{eqnarray*}
&&E\bigl\llvert h^{p}G_{n0}(Z_{t},Z_{0})
\bigr\rrvert^{2}\\
&&\qquad=16E\bigl\llvert h^{2p}\varepsilon
_{t}^{2}\varepsilon_{0}^{2}E_{0}^{2}
\bigl[ \sigma^{2}(%
\bar{X}_{0})W_{h}(
\bar{X}_{0},X_{t})W_{h}(\bar{X}_{0},X_{0})
\bigr] \bigr\rrvert
\\
&&\qquad\leq16h^{2p}E\bigl\llvert\varepsilon_{t}^{2}
\varepsilon_{0}^{2} \bigl[ E_{0}\sigma
^{2c}(\bar{X}_{0}) \bigr] ^{{2}/{c}} \bigl[
E_{0}W_{h}^{c}(%
\bar{X}_{0},X_{t})W_{h}^{c}(
\bar{X}_{0},X_{0}) \bigr] ^{{2}/{c}}\bigr\rrvert
\\
&&\qquad\leq16 h^{2p}C \bigl[ E\bigl\llvert\varepsilon_{t}^{4}
\varepsilon_{0}^{4}\bigr\rrvert\bigr] ^{{1}/{2}}
\bigl\{ E\bigl\llvert E_{0} \bigl( W_{h}^{c}(
\bar{X}_{0},X_{t})W_{h}^{c}(
\bar{X}_{0},X_{0}) \bigr) \bigr\rrvert^{{4}/{c}} \bigr
\} ^{{1}/{2}}
\\
&&\qquad\leq Ch^{2p} \bigl\{ E\bigl\llvert h^{-2cp+p}
\mathbf{A}_{c,h} ( X_{t},X_{0} ) \bigr\rrvert
^{{4}/{c}} \bigr\} ^{{1}/{2}}
\\
&&\qquad =O \bigl( h^{2p} \bigl[ h^{(-2cp+p)({4}/{c})z}h^{p} \bigr]
^{{1}/{2}%
} \bigr) =O \bigl( h^{({2}/{c}-{3}/{2})p} \bigr)
\end{eqnarray*}
for any $c>1$, where%
\begin{eqnarray*}
E_{0}\bigl[W_{h}^{c}(\bar{X}_{0},X_{t})W_{h}^{c}(
\bar{X}_{0},X_{0})\bigr] &=&h^{-2cp+p}\int
W_{h}^{c}(\bar{x}_{0},X_{t})W_{h}^{c}(
\bar{x}_{0},X_{0})f(%
\bar{x}_{0})\,d
\bar{x}_{0}
\\
&=&h^{-2cp+p}\mathbf{A}_{c,h}(X_{t},X_{0})
\end{eqnarray*}
by change of variable, where $\mathbf{A}_{c,h}(X_{t},X_{0})$ is a function
similar to $\mathbf{K}_{h}(X_{t}-X_{0})$. Thus, we obtain $%
\Vert h^{p}G_{n0}(Z_{t},Z_{0})\Vert_{2}\leq Ch^{({1}/{c}-
{3}/{4})p}$. By a
similar argument, we obtain the same order of magnitude for $%
\Vert h^{p}G_{n0}(Z_{t},\bar{Z}_{0})\Vert_{2}$. Thus, we have
$v_{n}(r)\leq Ch^{(%
{1}/{c}-{3}/{4})p}$, and condition (ii) holds, that is, $%
v_{n}(2)=o ( 1 )$, with $1<c<\frac{4}{3}$.

Next, to verify condition (iii), we shall evaluate $\Vert
h^{p}G_{n0}(Z_{0},\bar{%
Z}_{0})\Vert_{r}$ for $r<2+\frac{\delta_{0}}{4}$. By the Cauchy--Schwarz
inequality and change of variables, we have%
\begin{eqnarray*}
E\bigl\llvert h^{p}G_{n0}(Z_{0},Z_{0})
\bigr\rrvert^{r} &=&4^{\gamma}E\bigl\llvert h^{rp}
\varepsilon_{0}^{r}\varepsilon_{0}^{r}E_{0}^{r}
\bigl[ \sigma^{2}(\bar{X}_{0})W_{h}(
\bar{X}_{0},X_{0})W_{h}(\bar{X}_{0},X_{0})
\bigr] \bigr\rrvert
\\
&\leq&4^{\gamma}h^{2p}E\bigl\llvert\varepsilon
_{0}^{2r}\sigma^{2c}(\bar{X%
}_{0})^{{r}/{c}} \bigl[ E_{0}W_{h}^{2c}(
\bar{X}_{0},X_{0}) \bigr] ^{{%
r}/{c}}\bigr\rrvert
\\
&\leq&Ch^{rp}\bigl(E\varepsilon_{0}^{4r}
\bigr)^{{1}/{2}} \bigl[ E\bigl\llvert E_{0}W_{h}^{2c}(
\bar{X}_{0},X_{0})\bigr\rrvert^{{2r}/{c}} \bigr]
^{%
{1}/{2}}
\\
&=&O \bigl( h^{rp} \bigl[ h^{(1-2c)p\cdot{2r}/{c}} \bigr] ^{
{1}/{2}%
}
\bigr) =O \bigl( h^{rp({1}/{c}-1)} \bigr),
\end{eqnarray*}
where $E_{0}[W_{h}^{2}(\bar{X}_{0},X_{0})]=\int W_{h}^{2c}(\bar{x}%
_{0},X_{0})f_{\bar{X}_{0}}(\bar{x}_{0})\,d\bar{x}_{0}=O (
h^{(1-2c)p} ) $ by change of variable. Thus, we obtain $%
\Vert h^{p}G_{n0}(Z_{0},Z_{0})\Vert_{r}\leq Ch^{p(
{1}/{c}-1)}=Cn^{\omega p(1-%
{1}/{c})}$ given $h\propto n^{-\omega}$. Thus condition (iii)
holds by
choosing $c$ sufficiently small subject to the constraint of $c>1$.

Finally, we verify condition (iv). We first consider the case with
$t=0$ and
$s\neq0$. We have, by the Cauchy--Schwarz inequality and change of
variables,%
\begin{eqnarray*}
E\bigl\llvert h^{p}G_{ns}(Z_{0},Z_{0})
\bigr\rrvert^{2} &=&16E\bigl\llvert h^{2p}\varepsilon
_{0}^{2}\varepsilon_{0}^{2}E_{0}^{2}
\bigl[\bar{\varepsilon}_{s}%
\bar{\varepsilon}_{0}W_{h}(
\bar{X}_{s},X_{0})W_{h}(\bar{X}_{0},X_{0})
\bigr]\bigr\rrvert
\\
&\leq&16h^{2p}E\bigl\llvert\varepsilon_{0}^{4}
\bigl( E_{0}\bar{\varepsilon}%
_{s}^{c}
\bar{\varepsilon}_{0}^{c} \bigr) ^{{2}/{c}} \bigl[
E_{0}W_{h}^{c}(%
\bar{X}_{s},X_{0})W_{h}^{c}(
\bar{X}_{0},X_{0}) \bigr] ^{{2}/{c}}\bigr\rrvert
\\
&\leq&16h^{2p}\bigl(E\llvert\varepsilon_{0}\rrvert
^{8}\bigr)^{{1}/{2}}%
\bigl[ E\bigl\llvert
E_{0}^{{4}/{c}}W_{h}^{c}(
\bar{X}_{s},X_{0})W_{h}^{c}(%
\bar{X}_{0},X_{0})\bigr\rrvert\bigr] ^{{1}/{2}}
\\
&=&O \bigl( h^{2p} \bigl[ h^{2(1-c)p\cdot({4}/{c})} \bigr] ^{
{1}/{2}%
}
\bigr) =O \bigl( h^{2({2}/{c}-1)p} \bigr),
\end{eqnarray*}
where $E_{0}[W_{h}^{c}(\bar{X}_{s},X_{0})W_{h}^{c}(\bar
{X}_{0},X_{0})]=\int
W_{h}^{c}(\bar{x},X_{0})W_{h}^{c}(\bar{x}_{0}, X_{0})f_{\bar{X}_{s}\bar
{X}%
_{0}}(\bar{x},\break\bar{x}_{0})\,d\bar{x}\,d\bar{x}_{0}=O ( h^{2(1-c)p} ) $
by change of variable. Thus, we have $\Vert
h^{p}G_{ns}(Z_{0},\break Z_{0})\Vert_{2}\leq%
[ Ch^{2({2}/{c}-1)p} ] ^{{1}/{2}}=Ch^{({2}/{c}-1)p}$,
and so $n^{\gamma_{1}}\Vert h^{p}G_{ns}(Z_{0},Z_{0})\Vert
_{2}=\break n^{(1-{2}/{c}%
)\omega p+\gamma_{1}}$ if $h=O(n^{-\omega})$. Therefore, we obtain $%
\Vert h^{p}G_{ns}(Z_{0},Z_{0})\Vert_{2}=O(n^{-\gamma_{1}})$ with
$\gamma
_{1}= ( \frac{c}{2}-1 ) \omega p$, if we choose $c$ small enough
for $1<c<2$. For the case with $t\neq0$ and $s\neq0$, by a similar
argument, we have $\Vert h^{p}G_{ns}(Z_{t},Z_{0})\Vert_{2}\leq [
Ch^{2({2%
}/{c}-1)p} ] ^{{1}/{2}}=O ( h^{({2}/{c}-1)p} )$. Thus,
condition (iv) holds with $\gamma_{1}= ( \frac{c}{2}-1 ) \omega p$,
provided we choose $c$ small enough with $1<c<2$. Since all
conditions~(i)--(iv) hold, we have $V_{q}^{-1/2}h^{p/2}\hat{H}_{q}\stackrel{d}{%
\longrightarrow}N(0,1)$ by Tenreiro's (\citeyear{Ten97}) central limit
theorem.
\end{pf*}

\begin{pf*}{Proof of Theorem~\ref{th3}\normalfont{(ii)}}
We shall now derive the
asymptotic distribution of $\lambda_{n}$ under $\mathbb{H}_{n}(a_{n})$.
From Lemmas~\ref{le3.5} and~\ref{le3.6} and Propositions~\ref{pr3.3} and~\ref{pr3.4} below, we have
under $%
\mathbb{H}_{n}(a_{n})$,%
\[
\biggl[ \lambda_{n}-h^{-p}\sigma^{-2}c(K)\int
\sigma^{2}(x)\,dx \biggr] \Big/%
\sqrt{2\sigma^{-4}d(K)
\int\sigma^{4}(x)\,dx}\stackrel{d} {\longrightarrow}%
N(\xi,1).
\]
\upqed\end{pf*}

\begin{lemma}\label{le3.5}
Under the conditions of Theorem~\ref{th3},
$\lambda_{n}=\frac{n}{2}\frac{\operatorname{SSR}_{0}-\operatorname{SSR}_{1}}{\operatorname{SSR}_{1}}+o_{p}(h^{-p/2})$
under $\mathbb{H}_{n}(a_{n})$ with $a_{n}=n^{-1/2}h^{-p/4}$.
\end{lemma}

\begin{lemma}\label{le3.6}
Under the conditions of Theorem~\ref{th3}, $\hat{\sigma}_{n}^{2}\equiv n^{-1}\operatorname{SSR}_{1}=\sigma
^{2}+O_{p}(n^{-1/2})$
 under $\mathbb{H}_{n}(a_{n})$ with $a_{n}=n^{-1/2}h^{-p/4}$.
 \end{lemma}

\begin{proposition}\label{pr3.3}
Let $\widetilde{\operatorname{SSR}}_{0}$
and $\widetilde{\operatorname{SSR}}_{1}$ be defined in the same
way as
$\operatorname{SSR}_{0}$ and $\operatorname{SSR}_{1}$, respectively, with $\{
\varepsilon_{t}\}_{t=1}^{n}$ replacing $\{ \hat{\varepsilon}%
_{t}\}_{t=1}^{n}$. Then under\vspace*{1pt} the conditions of Theorem~\ref{th3},
$\operatorname{SSR}_{0}-\operatorname{SSR}_{1}=\widetilde{\operatorname{SSR}}_{0}-\widetilde{\operatorname{SSR}}_{1}+h^{-p/2}E[\delta
^{2}(X_{t})]+o_{p}(h^{-p/2})$ under $\mathbb{H}_{n}(a_{n})$
with $a_{n}=n^{-1/2}h^{-p/4}$.
\end{proposition}

\begin{proposition}\label{pr3.4}
Under the conditions of Theorem
\ref{th3} and $\mathbb{H}_{n}(a_{n})$ with $a_{n}=n^{-1/2}h^{-p/4}$,
\begin{eqnarray*}
&&\biggl[ \frac{\widetilde{\operatorname{SSR}}_{0}-\widetilde{\operatorname{SSR}}_{1}}{2\sigma^{2}}%
-h^{-p}\sigma
^{-2}c(K)\int\sigma^{2}(x)\,dx \biggr] \Big/\sqrt{2\sigma
^{-4}d(K)\int\sigma^{4}(x)\,dx}\\
&&\qquad\stackrel{d} {\longrightarrow
}N(\xi,1).
\end{eqnarray*}
\end{proposition}

\begin{pf*}{Proof of Lemma~\ref{le3.5}} Given in the supplementary material.
\end{pf*}

\begin{pf*}{Proof of Lemma~\ref{le3.6}} Given in the supplementary material.
\end{pf*}

\begin{pf*}{Proof of Proposition~\ref{pr3.3}} Given in the supplementary
material.~%
\end{pf*}

\begin{pf*}{Proof of Proposition~\ref{pr3.4}}
Proposition~\ref{pr3.4} follows from
Lemmas~\ref{le3.7} and~\ref{le3.8} below.
\end{pf*}

\begin{lemma}\label{le3.7}
Put $\hat{H}_{\lambda
}=n^{-1}\sum_{t=2}^{n}\sum_{s=1}^{t-1}H_{n}(Z_{t},Z_{s})$, where $%
Z_{t}=(\varepsilon_{t},X_{t}^{\prime})^{\prime}$, $%
H_{n}(Z_{t},Z_{s})=\varepsilon_{t}\varepsilon_{s}W_{h}(X_{t},X_{s})$
and
\[
W_{h}(X_{t},X_{s})= \biggl[
\frac{1}{f(X_{t})}+\frac{1}{f(X_{s})} \biggr] \mathbf{K}_{h}(X_{t}-X_{s})-
\int\frac{\mathbf{K}_{h}(X_{t}-x)\mathbf{K}%
_{h}(X_{s}-x)}{f(x)}\,dx.
\]
Suppose Assumptions~\ref{assa1} and~\ref{assa4} hold, $h\propto n^{-\omega}$
for $\omega\in(0,1/2p)$ and $p<4$.
Then
\[
\widetilde{\operatorname{SSR}}_{0}-\widetilde{\operatorname{SSR}}_{1}=h^{-p}
\biggl[ 2\mathbf{K}(0)-\int\mathbf{K}^{2}(\mathbf{u})\,d\mathbf{u}
\biggr] \int\sigma^{2}(x)\,dx+\hat{H}%
_{\lambda}+o_{p}
\bigl(h^{-p/2}\bigr).
\]
\end{lemma}

\begin{lemma}\label{le3.8}
Suppose Assumptions~\ref{assa1} and~\ref{assa4} hold,
and $h\propto n^{-\omega}$ for $\omega\in(0,1/2p)$. Define
\[
V_{\lambda}=2\int\biggl[ \mathbf{K}(\mathbf{u})-\frac{1}{2}\int
\mathbf{K}(%
\mathbf{v})\mathbf{K}(\mathbf{u}+\mathbf{v})\,d\mathbf{v}
\biggr] ^{2}\,du\int\sigma^{4}(x)\,dx.
\]
Then $V_{\lambda}^{-1/2}h^{p/2}\hat{H}_{\lambda}\stackrel{d}{
\longrightarrow}N(\xi,1)$.
\end{lemma}

\begin{pf*}{Proof of Lemma~\ref{le3.7}} Given in the supplementary material.
\end{pf*}

\begin{pf*}{Proof of Lemma~\ref{le3.8}} Given in the supplementary material.
\end{pf*}

\begin{pf*}{Proof of Theorem~\ref{th4}}
Pitman's asymptotic relative
efficiency of the $q_{n}$ test over the $\lambda_{n}$ test is the
limit of
the ratio of the sample sizes required by the two tests to have the same
asymptotic power at the same significance level, under the same local
alternative; see \citet{Pit79}, Chapter~7. Supposed $n_{1}$ and
$n_{2}$ are the
sample sizes required for the $q_{n}$ and $\lambda_{n}$ tests, respectively.
Then Pitman's asymptotic relative efficiency of $q_{n}$ to $\lambda
_{n}$ is
defined as%
%
\begin{equation}
\operatorname{ARE}(q_{n}\dvtx\lambda_{n})=\lim_{n_{1},n_{2}\rightarrow\infty}
\frac{n_{1}}{%
n_{2}} \label{eqA.5}
\end{equation}
under the condition that $\lambda_{n}$ and $q_{n}$ have the same asymptotic
power under the same local alternatives $n_{1}^{-
{1}/{2}}h_{1}^{-{p%
}/{4}}\delta_{1}(x)\sim n_{2}^{{1}/{2}}h_{2}^{-{p}/{4}}\delta
_{2}(x) $ in the sense that%
\[
\lim_{n_{1},n_{2}\rightarrow\infty}\frac{n_{1}^{-
{1}/{2}}h_{1}^{-{%
p}/{4}}\delta_{1}(x)}{n_{2}^{{1}/{2}}h_{2}^{-{p}/{4}}\delta
_{2}(x)}%
=1.
\]
Given $h_{i}=cn_{i}^{-\omega},i=1,2$, we have $n_{1}^{-2\gamma
}E[\delta
_{1}^{2}(X_{t})]\sim n_{2}^{-2\gamma}E[\delta_{2}^{2}(X_{t})]$, where
$%
\gamma=\frac{2-\omega p}{4}$. Hence,

\begin{equation}
\lim_{n_{1},n_{2}\rightarrow\infty} \biggl( \frac{n_{1}}{n_{2}} \biggr)
^{2\gamma}=\frac{E[\delta_{2}^{2}(X_{t})]}{E[\delta
_{1}^{2}(X_{t})]}. \label{eqA.6}
\end{equation}

On the other hand, from Theorem~\ref{th3}(ii), we have%
\[
\frac{\gamma(K)\lambda_{n_{1}}-\mu_{n_{1}}}{\sqrt{2\mu
_{n_{1}}}}\stackrel{%
d} {\rightarrow}N(\xi,1),
\]
under $\mathbb{H}_{n_{1}}(a_{n_{1}})\dvtx g_{0}(X_{t})=g(X_{t},\theta
_{0})+n_{1}^{-{1}/{2}}h_{1}^{-{1}/{4}}\delta_{1}(X_{t})$,
where $%
\xi=E[\delta_{1}^{2}(X_{t})]/\break [2\sigma^{-2}\sqrt{2 d(K)\int\sigma
^{4}(x)\,dx%
}]$. Also, from Theorem~\ref{th3}(i), we have%
\[
\frac{q_{n_{2}}-\nu_{n_{2}}}{\sqrt{2\nu_{n_{2}}}}\stackrel{d}
{\rightarrow}%
N(\psi,1)
\]
under $\mathbb{H}_{n_{2}}(a_{n_{2}})\dvtx g_{0}(X_{t})=g(X_{t},\theta
_{0})+n_{2}^{-{1}/{2}}h_{2}^{-{1}/{4}}\delta_{2}(X_{t})$,
where $%
\psi=E[\delta_{2}^{2}(X_{t})]/\break \sigma^{-2}\sqrt{2b(K)\int\sigma
^{4}(x)\,dx}$. To have the same asymptotic power, the noncentrality parameters must be
equal; namely $\xi=\psi$, or%
%
\begin{equation}
\frac{E[\delta_{1}^{2}(X_{t})]}{2\sqrt{2d(K)\int\sigma^{4}(x)\,dx}}=\frac{%
E[\delta_{2}^{2}(X_{t})]}{\sqrt{2b(K)\int\sigma^{4}(x)\,dx}}. %
\label{eqA.7}
\end{equation}
Combining (\ref{eqA.5})--(\ref{eqA.7}) yields%
\begin{eqnarray*}
\operatorname{ARE}(q_{n},\lambda_{n}) &=& \biggl[ \frac{2\sqrt{d(K)}}{\sqrt{b(K)}}
\biggr] ^{{1}/{(2\gamma)}}= \biggl[ \frac{4d(K)}{b(K)} \biggr] ^{
{1}/{(4\gamma)}}
\\
&=& \biggl[ \frac{\int ( 2\mathbf{K}(u)-\int\mathbf{K}(u)\mathbf{K}%
(u+v)\,du ) ^{2}\,dv}{\int ( \int\mathbf{K}(u)\mathbf{K}%
(u+v)\,du ) ^{2}\,dv} \biggr] ^{{1}/{(2-\omega p)}}.
\end{eqnarray*}
Finally, we show $\operatorname{ARE}(q_{n}\dvtx\lambda_{n})\geq1$ for any positive kernels
with $K(\cdot)\leq1$. For this purpose, it suffices to show
\[
\int\biggl[ 2\mathbf{K}(u)-\int\mathbf{K}(u)\mathbf{K}(u+v)\,du \biggr]
^{2}\,dv\geq\int\biggl[ \int\mathbf{K}(u)\mathbf{K}(u+v)\,du \biggr]
^{2}\,dv
\]
or equivalently,%
\[
\int\mathbf{K}^{2}(v)\,dv\geq\iint\mathbf{K}(u)\mathbf{K}(v)
\mathbf{K}%
(u+v)\,du\,dv.
\]
This last inequality follows from Zhang and Dette [(\citeyear{ZhaDet04}),
Lemma 2]. This
completes the proof.%
%
\end{pf*}
\end{appendix}

\section*{Acknowledgments}
We would like to thank the Editor, an Associate
Editor and two anonymous referees for insightful comments and suggestions
which significantly improved our paper. We also thank participants of
2008 International Symposium on Recent Developments in Time Series
Econometrics, Xiamen, China and 2008 Symposium on Econometric Theory and
Applications (SETA), Seoul, Korea, 2009 Time Series Conference, CIREQ,
Montreal, Canada, Midwest Econometrics Group Meeting, Purdue, and seminar
participants at Cornell University, Hong Kong University of Science and
Technology, Korea University, Ohio State University, Texas A\&M and UC
Riverside for their helpful comments and discussions.

\begin{supplement}[id=suppA]
\stitle{Supplementary material for a loss function approach
to model specification testing and
its relative efficiency}
\slink[doi]{10.1214/13-AOS1099SUPP} 
\sdatatype{.pdf}
\sfilename{aos1099\_supp.pdf}
\sdescription{In this supplement, we present the detailed proofs of
Theorems~\ref{th1}--\ref{th4} and
report the simulation results with the bandwidth $h=S_{X}n^{-1/5}$.}
\end{supplement}

\printaddresses

\end{document}